\documentclass[a4paper,12pt,leqno,oneside]{amsart}
\usepackage[latin1]{inputenc}
\usepackage[T1]{fontenc}

\usepackage{amssymb,calrsfs,ae}

\swapnumbers
\let\dsp=\displaystyle
\theoremstyle{plain}
\newtheorem{prop}{Proposition}[section]
\newtheorem{theo}[prop]{Theorem}
\newtheorem{cor}[prop]{Corollary}
\newtheorem{lem}[prop]{Lemma}
\newtheorem{conj}[prop]{Conjecture}

\theoremstyle{definition}
\newtheorem{defi}[prop]{Definition}

\theoremstyle{remark}

\newtheorem{rem}[prop]{Remark}

\newcommand{\thismonth}{\ifcase\month\or
  January\or February\or March\or April\or May\or June\or July\or
  August\or September\or October\or November\or December\fi
  \space\number\year}
\DeclareMathOperator{\Tr}{Tr}
\DeclareMathOperator{\im}{im}

\DeclareMathOperator{\spec}{spec}

\def\R{\mathbb R}
\def\C{\mathbb C}
\def\N{\mathbb N}
\def\Z{\mathbb Z}

\renewcommand{\epsilon}{\varepsilon}
\renewcommand{\Re}{\operatorname{Re}}
\renewcommand{\Im}{\operatorname{Im}}

\newcommand{\NM}{{\mathbb N}}
\newcommand{\RM}{{\mathbb R}}\newcommand{\setR}{\RM}

\newcommand{\ZM}{{\mathbb Z}}\newcommand{\setZ}{\ZM}

\newcommand{\vol}{\operatorname{vol}}
\newcommand{\Scal}{\operatorname{Scal}}

\newcommand{\tr}{\operatorname{Tr}}

\newcommand{\End}{\operatorname{End}}

\renewcommand{\geq}{\geqslant}
\renewcommand{\leq}{\leqslant}
\newcommand{\cerc}{{\mathbb S}}

\renewcommand{\exp}{\operatorname{e}}

\newcounter{mnotecount}[section]
\renewcommand{\themnotecount}{\thesection.\arabic{mnotecount}}
\newcommand{\mnote}[1]
{\protect{\stepcounter{mnotecount}}$^{\mbox{\footnotesize  $
      \bullet$\themnotecount}}$ \marginpar{\raggedright\tiny\em
    $\!\!\!\!\!\!\,\bullet$\themnotecount: #1} }

\begin{document}

\title[Diabatic limit and CR 3-manifolds]{Diabatic limit, eta
  invariants and Cauchy-Riemann manifolds of dimension $3$}

\author{Olivier Biquard}
\address{Institut de Recherche Mathématique Avancée\\ 
  CNRS et Université Louis Pasteur \\ 7 rue René Descartes, 67084
  Strasbourg Cedex\\ France}
\email{olivier.biquard@math.u-strasbg.fr} 

\author{Marc Herzlich}
\address{Institut de Mathématiques et de Modélisation de Montpellier\\
  CNRS et Université Mont­pellier~II\\ Place Eugène Bataillon, 34095 Montpellier Cedex 5\\ France}
\email{herzlich@math.univ-montp2.fr} 

\author{Michel Rumin}
\address{Laboratoire de Mathématiques d'Orsay\\ CNRS et Université Paris Sud\\ 
 91405 Orsay Cedex\\ France}
\email{michel.rumin@math.u-psud.fr} 


\begin{abstract}
  We relate a recently introduced non-local geometric invariant of compact
  strictly pseudoconvex Cauchy-Riemann (CR) manifolds of dimension $3$
  to various $\eta$-invar-iants in CR geometry: on the one hand a
  renormalized $\eta$-invariant appearing when considering a sequence
  of metrics converging to the CR structure by expanding the size of
  the Reeb field; on the other hand the $\eta$-invariant of the middle
  degree operator of the contact complex. We then provide explicit
  computations for a class of examples: transverse circle invariant CR
  structures on Seifert manifolds.  Applications are given to the
  problem of filling a CR manifold by a complex hyperbolic manifold,
  and more generally by a Kähler-Einstein or an Einstein metric.
\end{abstract}
\keywords{CR manifolds of dimension $3$, pseudohermitian structures, 
eta invariants, contact complex} 
\subjclass{32V05, 32V20, 53C20, 58J28}

\thanks{The second author is supported in part by a Young Researchers 
{\sc aci} program of the French Ministry of Research.}  

\maketitle

\section{Introduction}

In \cite{ob-mh1} the first two authors of this paper introduced a new
invariant, called the $\nu$-invariant, of strictly pseudoconvex
Cauchy-Riemann (CR) compact 3-manifolds. This invariant was obtained
by taking the limit of the $\eta$-invariants of an adequately defined
(but quite complicated) sequence of Riemannian metrics approximating
the CR structure, after cancellation of the possibly diverging terms
by adding well-chosen local contributions.  We claimed in
\cite{ob-mh1} that this invariant was an analogue in CR geometry of
the $\eta$-invariant in conformal geometry.  However, its rather
abstract definition makes it difficult to compute explicit expressions
for it or to get a further understanding of its properties. The
goal of this paper is then to provide links between $\nu$ and other
natural $\eta$-invariants in CR geometry.

\medskip

In a first step, we introduce a renormalized $\eta$-invariant that
takes into account the fact that CR geometry can be seen as a limit of
a sequence of conformal structures that diverges outside the contact
distribution. If a compatible contact form $\theta$ is fixed on the CR
manifold $M$, one considers the family of metrics
\begin{equation}\label{eq:gammaeps}
h_\epsilon = \epsilon^{-1} \theta^2 + \gamma , 
\end{equation}
where $\gamma=d\theta(\cdot,J\cdot)$ and $J$ is the underlying complex
structure on the contact distribution. \emph{When $\epsilon$ goes to
  $0$} the metrics $h_\epsilon$ blow up except in the contact
distribution, and therefore diverge to the Carnot-Carathéodory
metric associated to the CR structure and the contact form (this is
one of the main motivation for considering this kind of sequences). A
natural object one can consider is the constant term $\eta_0$ in an 
asymptotic expansion for $(\eta(h_\epsilon))$ in powers of
$\epsilon$, \emph{when $\epsilon$ goes to $0$}. This always exists, as
we shall see, and we shall call it the \emph{renormalized
  $\eta$-invariant of the pseudohermitian manifold $(M,\theta)$}. This
invariant is of course much more easily studied than the $\nu$-invariant, 
because it is 
built from the sequence \eqref{eq:gammaeps} of metrics that is much simpler
than the one used to build $\nu$ in \cite{ob-mh1}. Note however that
it is a pseudohermitian invariant, \emph{i.e.} it depends on the
choice of $\theta$, contrarily to $\nu$.

In the other direction, \emph{i.e. when $\epsilon$ goes to $\infty$}, on can also
obtain another natural invariant in case the Tanaka-Webster torsion of
$(M,\theta)$ vanishes, that is when the action of the Reeb vector field is
isometric. In this case, $\eta(h_\epsilon)$ converges and its limit
$\eta_{\mathrm{ad}}$ is the so called \emph{adiabatic limit}. It has
attracted much attention in the past few years, see
\cite{bismut-cheeger,dai} for instance.  We shall call the reverse
process of taking a limit when $\epsilon$ goes to $0$ a \emph{diabatic}
limit. When torsion vanishes, it turns out that the diabatic $\eta_0$
equals the adiabatic $\eta_{\mathrm{ad}}$.

Our first result shows that the difference between the CR invariant
$\nu$ and the pseudo-hermitian $\eta_0$ is an integral of a local
contribution involving the square of the Tanaka-Webster curvature.

\begin{theo}\label{eta0}
  For any compact strictly pseudoconvex Cauchy-Riemann $3$-manifold
  $M$, and any choice $\theta$ of contact form, one has
  \begin{equation}\label{eq:eta0} 
  \nu(M) \ = \ - \, 3 \,\eta_0(M,\theta) \, + \,\frac{1}{16\pi²}
  \,\int_{M} R^2 \theta\land d\theta \ ,
  \end{equation}
  where $R$ is the Tanaka-Webster curvature of $(M,\theta)$.
\end{theo}

This yields a new definition of the $\nu$-invariant, see Remark
\ref{rem:ndnu}, together with some explicit computations: they can be
done on manifolds on which $\eta_0$ is computable. We are then able to
apply this to transverse $\cerc^1$-invariant CR structures on Seifert
manifolds. The CR manifolds we are interested in come with a locally
free action of $\cerc^1$ that is transverse to the contact
distribution, and preserves both the contact and the complex
structures. We shall call them \emph{Cauchy-Riemann-Seifert manifolds}
(in short CR-Seifert). We refer to \cite{ka-tsu} for more information
on the more general class of $\cerc^1$-invariant CR structures.
CR-Seifert manifolds can be efficiently described as orbifold
$\cerc^1$-bundles over $2$-dimensional orbifolds. At each orbifold
point on the base, the orbifold bundle data consists of the following:
if the local fundamental group is $\setZ/\alpha\setZ$
($\alpha\in\NM^*$), a generator acts on a local chart around $p$ on
the basis manifold as $\exp^{i\frac{2\pi}{\alpha}}$ and on the
fiber as $\exp^{i\frac{2\pi\beta}{\alpha}}$ with $\beta$
prime to $\alpha$.  The orbifold $\cerc^1$-bundles are topologically
classified by their degrees (first Chern numbers), which are in this
case rational numbers. One then endows the manifold with an invariant
strictly pseudoconvex CR structure as follows: the underlying contact
structure is provided by an equivariant connection $1$-form on the
bundle, whereas the complex structure is induced from the basis
(orbifold) Riemann surface; the strict pseudoconvexity condition
constrains the degrees $d$ of these $\cerc^1$-bundles to be negative.

\begin{theo}\label{maintheo1}
  Let $M$ be a compact strictly pseudoconvex CR-Seifert $3$-manifold,
  of degree $d$ over the orbifold surface $\Sigma$, and with
  $\cerc^1$-action generated by the Reeb field of a contact form
  $\theta$. If $R$ is the Tanaka-Webster curvature of $(M,\theta)$,
  then
 \begin{equation}\label{eq:nu} 
  \nu(M) = - d - 3  - 12 \sum_{j=1}^{p} s(\alpha_j,1,\beta_j)
  + \frac{1}{8\pi} \int_\Sigma R^2  d\theta  \ ,
  \end{equation}
  where $s(\alpha,\rho,\beta)$ is the Rademacher-Dedekind sum
  $\frac{1}{4\alpha} \sum\limits_{k=1}^{\alpha -1} \cot \left(
    \frac{k\rho\pi}{\alpha}\right) \cot \left(
    \frac{k\beta\pi}{\alpha}\right)$ .
\end{theo}

The Tanaka-Webster curvature $R$ of such an $(M,\theta)$ actually
coincides with Riemannian curvature of the base $\Sigma$, if it is
endowed with the metric $\gamma = d\theta ( \cdot, J \cdot)$.  When
this curvature is constant, \eqref{eq:nu} specializes into the
following interesting formula, which shows that the $\nu$-invariant is
a topological invariant in this case:

\begin{cor}
  \label{maintheo3}
  Let $M$ be a CR-Seifert manifold as above, with constant
  Tanaka-Webster curvature. Let $\chi$ be the rational Euler
  characteristic of $\Sigma$. Then,
  \begin{equation}
    \label{eq:nuconstant} 
    \nu(M) = - d - 3  - \frac{\chi^2}{4d} 
    - 12 \sum_{j=1}^{p} s(\alpha_j,1,\beta_j) \ .
  \end{equation}
\end{cor}

\medskip

However, Theorem \ref{eta0} is not entirely satisfactory, as it
provides a link between the CR invariant $\nu$ and the diabatic
invariant $\eta_0$; one would instead prefer a relationship between $\nu$
and invariants defined directly in terms of the CR or pseudohermitian
geometry. One such object is the contact-de Rham complex introduced in
\cite{Rumin94}, and especially the $\eta$-invariant of the middle degree
operator appearing there.

The relevant operator (denoted by $D*$ henceforth) is the analogue in
this setting of the {\sl boundary operator for the signature} $±(d* -
*d)$ that gives rise to the $\eta$-invariant on $3$-dimensional
Riemannian manifolds.  It is known that the spectrum of the operator
$D*$ appears in the rescaled limit of the collapsing spectrum of
$P_\epsilon= ±(d*_\epsilon - *_\epsilon d)$ for the metrics
$h_\epsilon$ of \eqref{eq:gammaeps}, when performing the diabatic
limit \cite{Rumin00}. However, this limit is not uniform enough to
yield a direct relation between the $\eta$-invariants. In this paper,
we prove a general relation between $\nu$ and $\eta(D*)$ in the
special case provided by our transverse $\cerc^1$-invariant CR
manifolds.  In effect, we show that $\eta(D*)$ and $\nu$ differ only
by a simple local term in the Tanaka-Webster curvature of any chosen
pseudohermitian structure. Our second main set of results then reads:

\begin{theo}\label{th:D*eta0}
  Let $M$ be a compact strictly pseudoconvex CR-Seifert $3$-manifold,
  with $\cerc^1$-action generated by the Reeb field of a contact form
  $\theta$. If $R$ is the Tanaka-Webster curvature of $(M,\theta)$ and
  $D$ is the middle operator of the contact complex, then
\begin{equation}\label{eq:D*eta0}
\eta_0(M,\theta) = \eta(D*) + \frac{1}{512} \,\int_M  R^2 \theta\land d\theta  .
\end{equation} 
\end{theo}

\begin{cor}\label{maintheo2}
  Let $M$ be a CR-Seifert $3$-manifold as above, then one has:
  \begin{equation}\label{eq:etanu} 
  \nu(M) \ = \ - \,3 \,\eta(D*) \  + 
  \ \left(\frac{1}{16\pi²} - \frac{3}{512}\right)
  \,\int_{M} R^2 \,\theta\land d\theta \ .
  \end{equation}
\end{cor}

The philosophy underlying our results is indeed the following: whereas
$\nu$ is easily related to $\eta_0$, $\eta(D*)$ compares itself more
easily with $\eta_0$ rather than to $\nu$.  This somehow ``explains''
the quite strange combination of constants appearing in front of the
curvature term in \eqref{eq:etanu} in theorem \ref{maintheo2}: it is a
sum of diabatic contribution steming from theorem \ref{eta0} and a
purely spectral term linking $\eta(D*)$ and $\eta_0$, as will be
apparent from section \ref{sec:metrics-h_0r-from}.

\medskip

For general CR manifolds, we expect that when we take the diabatic
limit $\varepsilon\to0$, the collapsing spectrum of $P_\varepsilon$ gives the contribution
$\eta(D*)$ in the limit, while the remaining part of the spectrum, after
renormalization, gives only an integral of local terms. This leads to
the following conjecture.

\begin{conj}\label{mainconj}
  There exists a constant $C$ such that, for any compact strictly
  pseudoconvex Cauchy-Riemann $3$-manifold $M$ and any choice $\theta$
  of contact form, one has
  \begin{equation}\label{eq:etanugeneral} 
  \nu(M) \ = \ -\, 3 \,\eta(D*) \ +
  \ \left(\frac{1}{16\pi²} - \frac{3}{512}\right)
  \,\int_{M} R^2 \theta\land d\theta \ + \ C \,\int_{M} |\tau|^2
  \,\theta\land  
  d\theta \ ,
  \end{equation}
  with $R$ and $\tau$ the Tanaka-Webster curvature and torsion of
  $(M,\theta)$.
\end{conj}

As a first indication for the conjecture, we shall give in Theorem 
\ref{maintheo4} an abstract argument that shows that there exists a 
CR invariant of the form $\eta(D*)+C_1\int R^2 + C_2\int |\tau|^2$.
Unfortunately, we are unable to calculate the constants completely,
see Remark \ref{maintheo5}.

It is known that the $\eta$-invariant of the boundary operator for
signature is conformally invariant. If the conjecture is true, then
this is no more the case for $\eta(D*)$, which is a priori an
invariant of the pseudo-hermitian structure only: it depends on the
choice of a metric in the conformal class adapted to the CR structure.

\medskip
 
A third goal of this paper is to provide some geometric applications
on CR-Seifert manifolds, mainly with constant curvature. They are
spherical (locally isomorphic to the standard CR sphere $\cerc^3$),
hence are the boundary at infinity of a complex hyperbolic metric
defined in a neighbourhood $(0,\epsilon)× M$ of $M$ (in the case of
the $3$-sphere we can of course extend the metric globally to get the
Bergmann metric on the 4-ball). From \cite[Theorem 1.2]{ob-mh1} and
Theorem \ref{maintheo3}, we get the following obstruction for this
neighbourhood to have a global extension to a smooth complex
hyperbolic surface (with only one end):

\begin{cor}\label{cor:1}
  If a CR-Seifert manifold $M^3$ is the boundary at infinity of a
  complex hyperbolic metric defined on the interior of a smooth
  compact manifold $N^4$ with boundary $M$, then one has necessarily $
  \nu(M) = - \chi(N) + 3 \tau(N)$, where $\chi(N)$ and $\tau(N)$
  denote the Euler characteristic and signature of $N$.  In
  particular, $\nu(M)$, as provided by the formula
  \emph{(\ref{eq:nu})}, is an integer.
\end{cor}

This is a topological constraint on a filling, which we can restate in
the smooth case (no orbifold singularities):

\begin{cor}\label{cor:2}
  Let $M$ be a $\cerc^1$-bundle of degree $d$ over a Riemann surface
  $\Sigma$ of Euler characteristic $\chi$, with a $\cerc^1$-invariant
  spherical CR structure.  If $ \frac{ \chi^2}{4d}$ is not an integer
  then $M$ is not the boundary at infinity of a complex hyperbolic
  metric.
\end{cor}

The case $d=\frac{\chi}{2}$ yields an integer, and indeed, if $\Sigma$
is hyperbolic, $N$ can be taken to be the disk bundle of a square root
of the tangent bundle of $\Sigma$, which is well known to carry a
complex hyperbolic metric issued from a representation of
$\pi_1(\Sigma)$ in $SU(1,1) \subset SU(1,2)$. Our obstruction then
gives an interesting hint on whether a spherical CR-Seifert
$3$-manifold may appear as the quotient of the complement of the limit
set in the $3$-sphere of some discrete fixed point-free subgroup of
$SU(1,2)$ \cite{apanasov}.

More generally, the calculation in Theorem \ref{maintheo1} gives an
obstruction for $M$ to be the boundary at infinity of a
Kähler-Einstein or Einstein metric. The manifolds considered in this
paper are known to bound a complex Stein space with at most a finite
number of singular points \cite{harvey-lawson} and one may wish to
endow it with a Kähler-Einstein metric as in Cheng-Yau
\cite{cheng-yau}.  The type of metric to be considered has the same
kind of asymptotic expansion near the boundary $M$ as the Bergman
metric \cite{biquard-symmetric}; we called them ``asymptotically
complex hyperbolic'' (ACH) in \cite{ob-mh1}. If no singular points are
present and if the Cheng-Yau metric exists, one gets from the
Miyaoka-Yau inequality proved in \cite{rollin} the following:

\begin{cor}\label{cor:3}
  Let $M$ be as in Theorem \ref{maintheo1}. If $M$ is the boundary at
  infinity of an ACH Einstein metric on $M^4$, such that a
  Kronheimer-Mrowka invariant of $(N,M)$ is nonzero \upn{(}in
  particular, if $M$ is the boundary at infinity of a Kähler-Einstein
  metric on $N$\upn{)}, then
  $$
  \chi(N) - 3 \tau(N) \geqslant -\nu(M) = d + 3 + 12 \sum_{j=1}^{p}
  s(\alpha_j,1,\beta_j) - \frac{1}{8\pi}\,\int_\Sigma R^2 \,
  d\theta \ .$$
\end{cor}

For more information on Stein fillings, see
\cite{lisca-matic,stipsicz1}. The Kronheimer-Mrowka invariants are
Seiberg-Witten type invariants defined for a compact 4-manifold with
contact boundary; in particular, they do not vanish if $M$ carries a
symplectic form compatible with the contact structure on the boundary,
and this implies the Miyaoka-Yau inequality \cite{rollin}. This
inequality can of course be obtained directly for Kähler-Einstein
metrics.

\medskip

The paper is organized as follows. After recalling the definition of
the $\nu$-invariant in section \ref{sec:nu}, we define the
renormalized $\eta$-invariant $\eta_0$ and compare it to $\nu$ in
sections \ref{sec:adia} and \ref{sec:2:3}. The proof relies on
relatively simple considerations on $\eta$-invariants and Chern-Simons
theory, that prove that the difference between $\nu + 3\,\eta_0$ is
necessarily of the expected form: an integral term in the square of
the curvature and the squared norm of the torsion. The constants in
front of these local terms are then computed by considering
sufficiently many examples: left invariant structures on the
$3$-sphere.

The reader will then find in section \ref{sec:Sfrt} the explicit
computations of $\nu$ on CR-Seifert manifolds.
  
Taking one step further, sections \ref{sec:4:1} to
\ref{sec:decompostion-sp-d} lead to the relation between $\eta_0$ and
$\eta(D*)$ in the case of transverse $\cerc^1$-invariant CR
structures. The proof of Theorem \ref{maintheo2} relies on an explicit
study of the spectra of the $D*$ operator and the boundary operator
for the signature $±(d_{\epsilon}* - *d_{\epsilon})$ on closed
$2$-forms for the sequence of Riemannian metrics $h_\epsilon$ that
performs the diabatic limit in \eqref{eq:gammaeps}.  This can be done
only for $\cerc^1$-invariant structures and index theory shows once
again that a relation of the expected type must exist. One then has
again to evaluate the constant in front of the integral term by
looking at explicit computations of both $\eta(D*)$ and $\nu$ on the
standard sphere.

The existence of a CR invariant of the form 
$\eta(D*)+C_1\int R^2 + C_2\int |\tau|^2$ is considered in section
\ref{sec:eta+}. We also present a proof of the existence of $\eta(D*)$
on any compact strictly pseudoconvex CR manifold of dimension $3$, a
fact certainly known to specialists but whose proof seems to have
never been published so far.

The paper ends with a short section \ref{sec:coro} devoted to the
proof of the corollaries and to some generalizations, and also to a
comparison with the results one can get in the Kähler-Einstein case
using the $\mu$-invariant of Burns and Epstein \cite{BE88}.

\bigskip

\section{The $\nu$-invariant}\label{sec:nu}

Let $M$ be a $3$-dimensional compact strictly pseudoconvex CR
manifold, \emph{i.e.} a compact manifold $M$ endowed with a complex
structure $J$ defined on a contact distribution $H$ in $TM$.

A pseudohermitian structure $(M,\theta)$ consists in the additional
choice of a contact form $\theta$. It induces a metric $\gamma =
d\theta(\cdot,J\cdot)$ on $H$ and a splitting of both $TM$ and
$T^*M$ by means of the Reeb vector field $T$ defined by $\theta(T)=1$
and $\iota_Td\theta=0$. The Tanaka-Webster connection is then defined
by working in a local coframe $(\theta,\theta^1,\theta^{\bar 1})$ such
that $d\theta=i\theta^1\land \theta^{\bar 1}$: the connection form is
a purely imaginary 1-form $\omega^1_1$, and the torsion $\tau^1$ is a
$(0,1)$-form such that
$$
d\theta^1 = \theta^1\land \omega^1_1 + \theta\land \tau^1 , $$
and
the curvature $R$ is defined by
$$
d\omega^1_1 = -iR\, d\theta + (\tau^{\bar
  1}_{,\bar{1}}-\tau^1_{,1})\land \theta . $$
In more invariant terms,
it is the only metric and complex compatible connection $\nabla^W$ on $H$ such
that the torsion $\tau = T^{\nabla^W}(T,\cdot)_{|H}$ anticommutes with $J$.

\medskip

Given a pseudohermitian manifold $(M,\theta)$, one can define a first
metric $g_0$ on the product $N = \setR_+ × M$ by
\begin{equation}
  \label{eq:2:1}
  g_0 = dr^2 + h_0(r) , \qquad h_0(r) = \exp^{2r} \theta^2 + \exp^r \gamma .
\end{equation}
Here we think of the initial $M$ as a boundary of $M$ at infinity
(\emph{i.e.} when $r$ goes to infinity). Remark that when doing a conformal change
$\theta'=f\theta$, one gets a metric $g_0'=(dr')^2+\exp^{2r'}
f^2\theta^2 + \exp^{r'} f\gamma$, and the difference $g_0'-g_0$ goes
to zero at infinity after the coordinate change $r=r'+\log f$.
Therefore the asymptotic behaviour of the metric $g_0$ depends only on
the CR structure. We note moreover that
$$
h_0(r) = \exp^r \, ( \exp^r \theta^2 + \gamma) = \epsilon^{-1}
h_\epsilon ,
$$
where $h_\epsilon$ is the metric introduced in equation
\eqref{eq:gammaeps}, with $\epsilon = \exp^{-r}$.

One can extend $J$, initially defined on $M$, to an almost complex
structure $J_0$ on $N$ by defining
$$
J_0 \partial_r = \exp^{-r} T,
$$
where $T$ is the Reeb field associated to $\theta$.  As explained
in \cite{ob-mh1} the curvature of $g_0$, together with $J_0$, is
asymptotic when $r \to +\infty$ to curvature of the complex hyperbolic
plane with holomorphic sectional curvature $-1$.

One can add higher order corrections to $J_0$ and $g_0$ to get a {\sl
  uniquely defined} jet of Kähler-Einstein metric $g_{KE}$ up to order
$\exp^{-2r}$ (relatively to $g_0$), when $r$ tends to infinity. This
development is expressed with the covariant derivatives of
Tanaka-Webster curvature and torsion $R$ and $\tau$ of the
pseudohermitian manifold $(M,\theta)$, and calculated in \cite[theorem
3.3 and corollary 3.4]{ob-mh1}. More precisely, one finds an infinite
series $J(r) =J_0+J_1\exp^{-r}+J_2\exp^{-2r}+\cdots$ giving an integrable
(formal) complex structure $J(r)$ on $N$, whose first terms are
\begin{displaymath}
  J (r) = J_0 - 2 \exp^{-r} \tau + \exp^{-2r} ( 2 |\tau|^2 - J_0
  \nabla_T\tau ) + \cdots , 
\end{displaymath}
and a unique \emph{finite} jet of Kähler-Einstein metric $g_{KE}$ on
$M$, that is locally determined up to order 2 by $(M,\theta)$: given
some choice of coframe $\theta^1 \in \Omega^{1,0}H$, the expression of
its Kähler form $\omega$ is
\begin{equation}\label{eq:ke}
    \begin{split}
      \omega & = \exp^r (dr \land \theta + d \theta) - \frac{R}{2} d
      \theta  \\
      & \quad + \frac{4}{3} \bigl(\frac{i}{8} R_{,\bar 1} \vartheta^0
      \land \theta^{\bar 1} - \frac{i}{8} R_{, 1} \vartheta^{\bar 0}
      \land \theta^1 - \frac{1}{2} \tau^1_{\bar 1, 1} \vartheta^0
      \land \theta^{\bar 1} - \frac{1}{2} \tau^{\bar 1}_{1, \bar 1}
      \vartheta^{\bar 0} \land \theta^1 \bigr) \\
      & \quad - \frac{\Delta_H R}{2} \exp^{-r} d\theta - \frac{2}{3}
      \Bigl( \frac{R^2}{8} - |\tau |^2 - \frac{\Delta_H R}{6} +
      \frac{2 i}{3} ( \tau^1_{\bar 1, 11} - \tau^{\bar
        1}_{1, \bar 1 \bar 1}) \Bigr) \exp^{-r} dr \land \theta \\
      & \quad + \frac{2}{3} \Bigl( \frac{R^2}{8} - |\tau |^2 -
      \frac{\Delta_H R}{12} - \frac{i}{3} ( \tau^1_{\bar 1, 11} -
      \tau^{\bar 1}_{1, \bar 1 \bar 1}) \Bigr) \exp^{-r} d \theta +
      o(\exp^{-2r}),
\end{split}
\end{equation}
where $(\vartheta^0 = \exp^{-r} dr + i \theta,\vartheta^1)$ is a
coframe of $\Omega^{1,0}N$ associated to $J(r)$.

It is explained in \cite{ob-mh1} why higher order terms in $\omega$
are irrelevant in all what concerns the $\nu$-invariant to be defined
below.  We will denote by $g_{KE}$ the metric on $N$ given by this
second order jet of Kähler metric $g_{KE} = \omega(\cdot,J(r)\cdot)$.  We
observe that $g_{KE}$ has an \emph{universal} polynomial expression in
the powers of $\exp^{r}$, with coefficients that are tensorial in the
covariant derivatives of $R$ and $\tau$. By construction the leading
term of $g_{KE}$ is $g_0$ as given in \eqref{eq:2:1}, and the family
of metrics $h(r)$ induced on
$$
M_r = \{r\}× M \simeq M
$$
is asymptotic to $h_0(r)$ in \eqref{eq:2:1}.

Finally, an important point here is that, although we have chosen a
contact form to write down the formulas for $g_{KE}$, actually it does
depend only on the CR structure, not on the pseudohermitian structure.
This is because the filling complex structure on $N$ depends only on
$J$, as does the zeroth order term of $g_0$, and the finite part of
the Kähler-Einstein metric that we need is uniquely determined.

\medskip

We can now define the $\nu$-invariant of $M$. According to
\cite{ob-mh1}, it is obtained by taking the limit as $r$ goes to
infinity ({\sl i.e.} by taking the {\sl diabatic limit}) of the
boundary contribution on $M_r$ of the Atiyah-Patodi-Singer formula for
the characteristic number $\chi - 3 \tau$ of $[r_0, r] × M \subset N$,
\emph{with respect to the metric} $g_{KE}$.

\begin{defi}
  \label{defn:nu}
  The $\nu$-invariant of $M$ is
  \begin{displaymath}
    \nu(M) = \lim_{r \to +\infty}  B(g_{KE}, M_r) - 3 \eta(h(r)) , 
  \end{displaymath}
  where $\eta(h(r))$ is the $\eta$-invariant of the boundary operator
  for the signature $ S = (-1)^p (*d -d *) $ on $\Omega^{2p} M_r$ (see
  \cite{aps1}) with the metric $h(r)$, and $B(g_{KE}, M_r)$ is an
  integral over $M_r$ of the relevant secondary characteristic class,
  tensorially constructed from the curvature of $g_{KE}$ and the
  second fundamental form of $M_r \subset N$.
\end{defi}

It is shown in \cite{ob-mh1} that this limit exists and actually gives
rise to a CR invariant of $M$ (independent on the choice of the
contact form $\theta$).  The interested reader is referred to
\cite[(7.7)]{ob-mh1} for the general formula. We will not need the
precise form of the correction term $B(g_{KE}, M_r)$ in this paper.

\bigskip

\section{The renormalized eta invariant}\label{sec:adia} 

From its very definition, the invariant $\nu$ is a renormalisation of
$\eta$-invariants of a jet $h(r)$ of the very natural Kähler metric
$g_{KE}$ restricted to slides of large radii.  However, these metrics
are quite intricate (as formula \eqref{eq:ke} obviously shows), and
$\nu$ itself is given by a limit of some complicated expressions built
from these metrics.  For these reasons we would like to describe how
$\nu$ is related to the $\eta$-invariants of the much simpler
contact-rescaling family of metrics of formula \eqref{eq:gammaeps}:
$$
h_\epsilon=\epsilon^{-1} \theta^2 + \gamma .
$$
This can be done as follows: although $\eta$ is a priori not
locally computable from the metric, its variation is. Indeed from the
Atiyah-Patodi-Singer formula \cite{aps1} and Chern-Simons' theory
\cite{Ch-Si} one has
\begin{equation}
  \label{eq:2:3}
  \eta (h_{\epsilon_1}) - \eta (h_{\epsilon_0}) =
  \frac{1}{3} \int_M Tp_1 (\nabla_{\epsilon_1} , \nabla_{\epsilon_0}),  
\end{equation}
where $T p_1(\nabla_{\epsilon_1} , \nabla_{\epsilon_0})$ is
Chern-Simons' transgression form of the first Pontrjagin class relative
to the Levi-Civita connections of the \emph{product} metrics
$$
\widetilde g_\epsilon = dr^2 + h_\epsilon \quad \text{on} \quad N=
\setR × M.
$$
If $\nabla_{\epsilon_1} = \nabla_{\epsilon_0} + \alpha$ and
$\Omega_t$ is the curvature $2$-form of $\nabla_{\epsilon_0} + t
\alpha$, then
\begin{equation}\label{eq:2:4}
  T p_1(\nabla_{\epsilon_1} , \nabla_{\epsilon_0}) \  = \ 2
  \int_0^1 P_1(\alpha, \Omega_t) dt \  = \ -\frac{1}{4\pi^2} 
  \int_0^1 \Tr ( \alpha \land \Omega_t) dt.
\end{equation}
This leads quickly to the following lemma.

\begin{lem}\label{lem:etahom}
  Let $(M^3,J,\theta)$ be a strictly pseudoconvex pseudohermitian
  manifold, with metric $\gamma=d\theta(\cdot,J\cdot)$ on the
  contact distribution. Then the $\eta$-invariants of the family of
  metrics $h_\epsilon=\epsilon^{-1}\theta^2+\gamma$ have a
  decomposition in homogeneous terms:
  \begin{equation}
   \eta(h_\epsilon) = \sum_{i=-2}^2 \eta_i(M,\theta) \epsilon^i
   . \label{eq:2:6} 
   \end{equation}
   The terms $\eta_i$ for $i\neq0$ are integral of local
   pseudohermitian invariants of $(M,\theta)$, and the $\eta_i$ for
   $i>0$ vanish when the torsion vanishes.
\end{lem}

\begin{proof}
  Denote by $\nabla^W$ the Tanaka-Webster connection, with $\tau$
  being the torsion seen as a trace-free symmetric endomorphism of
  $H=\ker \theta$, $\tau^1$ (resp. $\tau^{\bar 1}$) being its
  expression as a $(0,1)$-form (resp $(1,0)$-form) relative to a
  choice of complex coframe $\theta^1$. One computes easily the
  difference $a=\nabla_\epsilon-\nabla^W$ (see the formulas in
  \cite[page 316]{Rumin94}), and the result is a decomposition into
  homogeneous terms of degrees $-1$, $0$ and $1$:
\begin{equation}\label{2:6bis}
  \nabla_\epsilon -\nabla^W = a = \sum_{-1}^1 a^{(i)} \epsilon^i ,
\end{equation}
where each $a^{(i)}$ is locally defined by the pseudohermitian
structure: $a^{(0)}$ and $a^{(-1)}$ are horizontal, but $a^{(1)}$ is
vertical, more precisely, for horizontal $X,Y\in H$ one has
\begin{align*}
  a^{(1)}_XY &= - \gamma(\tau(X),Y) T,\\
  a^{(0)}_XT &= \tau(X), \\
  a^{(-1)}_TY &= \frac{1}{2}JY .
\end{align*}
The output is the following decomposition for the curvature
\begin{align}
  \Omega(\nabla_\epsilon) & = \Omega(\nabla^W)+d^Wa+a\land a \\
  & = \sum_{-1}^1 \Omega^{(i)} \epsilon^i .
\end{align}
Indeed, the terms $\Omega^{(±2)}=a^{(±1)}\land a^{(±1)}$ clearly
vanish.  Moreover,
$$
\Omega^{(1)}=da^{(1)} + a^{(1)}\land a^{(0)} + a^{(0)}\land a^{(1)}
$$
vanishes when the torsion vanishes.  From equation \eqref{eq:2:4}
one has
\begin{equation*}
  \epsilon\frac{d}{d\epsilon}\eta(h_\epsilon) \ = \ 
  -\frac{1}{12\pi^2} \, \int_M \Tr( \Omega \land
  \epsilon\frac{da}{d\epsilon})\ 
  = \ \sum_{\substack{-2\leq i\leq 2\\i\neq0}} i \eta_i \,\epsilon^i
\end{equation*} 
where the $\eta_i$ ($i\neq0$) are local pseudo-hermitian invariants.
When the torsion vanishes, $a^{(1)}$ and $\Omega^{(1)}$ vanish, so
that $\eta_i$ vanishes for each $i>0$.
\end{proof}

From the conformal invariance of the $\eta$-invariant, one deduces
moreover immediately that, for a real number $\lambda>0$,
\begin{equation}
  \label{eq:2:7}
  \eta_i(M, \lambda \theta) = \lambda^{-i}\, \eta_i(M, \theta).
\end{equation}
so that $\eta_0(M,\theta)$ is scale (but not conformally) invariant.

\medskip

\begin{defi}
  \label{defi:eta_0}
  Let $(M^3,\theta)$ be a compact strictly pseudoconvex
  pseudohermitian manifold. The \emph{renormalized $\eta$-invariant of
    $(M,\theta)$} is the constant term $\eta_0(M,\theta)$ in the
  expansion \eqref{eq:2:6} for the $\eta$-invariants of the family of
  metrics $h_\epsilon=\epsilon^{-1}\theta^2+d\theta(\cdot,J\cdot)$.
\end{defi}

In the case where the torsion of $(M,\theta)$ vanishes, the terms
$\eta_i(M,\theta)$ in \eqref{eq:2:6} for $i>0$ vanish, so that, when
$\epsilon$ goes to infinity instead of $0$, one has
\begin{equation}
  \label{eq:1}
  \eta_0(M,\theta) = \lim_{\epsilon \to \infty}
  \eta(h_\epsilon) : = \eta_{\mathrm{ad}}.
\end{equation}
This corresponds to the geometric situation when the Reeb flow
preserves the metric. Then, when $\epsilon\to\infty$, the family of
metrics $h_\epsilon$ collapses with bounded connection and curvature.
This is the well-known \emph{adiabatic limit}, and $\eta_0(M,\theta)$
is then the adiabatic limit $ \eta_{\mathrm{ad}}$ of the
$\eta$-invariant. It has been much studied, in particular in the
geometrically meaningful situation when the Riemannian flow comes from
some fibration in circles over a surface
\cite{bismut-cheeger,dai}.

However, we are more interested in this paper in the opposite
direction: the {\sl diabatic limit}, or equivalently the case where
$\epsilon$ goes to $0$.  Although we will not need its precise
expression, making the calculations in the proof of lemma
\ref{lem:etahom} explicit shows the term $\eta_{-2}(M, \theta)$ never
vanishes on contact manifolds, and has to be of type $C\, \int_M
\theta \land d\theta$ for some universal non-zero constant $C$.
Therefore $\eta(h_\epsilon)$ always diverges at speed $\epsilon^{-2}$
in the diabatic limit, but the constant term $\eta_0(M,\theta)$ is
still well-defined. We called it the \emph{renormalized} $\eta$-invariant, as
it is reminiscent of other similar contexts where renormalized
invariants have been defined \cite{graham-volume,mh-vol,seshadri-vol}.

\bigskip

\section{The relation between $\nu$ and $\eta_0$}
\label{sec:2:3}

Our goal now is to prove Theorem \ref{eta0}, \emph{i.e.} to show that
on any CR manifold the $\nu$-invariant is related to $\eta_0$ in a
simple way.

\begin{lem}
  \label{lem:2:3} 
  There exist two constants $C_1$ and $C_2$ such that for any CR
  strictly pseudoconvex pseudohermitian manifold $(M^3,J,\theta)$,
  one has
  \begin{equation}
    \label{eq:2:9}
    \nu(M) + 3\eta_0(M, \theta) = C_1 \int_M R^2 \,\theta\land d\theta 
    \ + \ C_2 \int_M |\tau|^2\,\theta\land d\theta,  
  \end{equation}
  where $\eta_0(M, \theta)$ is the renormalized $\eta$-invariant of
  $(M,\theta)$, and $R$, $\tau$ are the Tanaka-Webster curvature and
  torsion of $M$.
\end{lem}

One can therefore look at $-\nu(M) /3$ as a local CR-conformal
correction of $\eta_0(M,\theta)$ (recall that $\eta_0(M,\theta)$ is
{\sl a priori} only invariant under the rescaling $\theta \to
\lambda\theta$ for $\lambda$ \emph{constant}).

\smallskip

\begin{proof}
  The metrics $g_{KE}$ and $h(r)=g_{KE}|_{\{r\}×M}$ issued from
  \eqref{eq:ke} are quite complicated, but are corrections of the
  model metrics $g_0$ and $h_0(r)$ defined in \eqref{eq:2:1}. More
  precisely, their expressions are universal polynomials in $\exp^r$
  and pseudohermitian invariant of $(M,\theta)$, and they do not
  actually depend on the choice of framing (except $\theta$) and the
  constants in front of each such term are universal, {\sl i.e.}
  independent of the manifold.  Therefore, using a transgression
  formula as in \eqref{eq:2:3} and \eqref{eq:2:4}, but between $h(r)$
  and $h_0(r)$, we see that $\eta(h(r)) - \eta(h_0(r))$ has to be an
  invariant universal expression of type
  \begin{equation}
    \label{eq:2:10}
    \sum_{k= -n}^n \exp^{kr} \int_M P_k(R, \tau, \nabla R , \nabla \tau,
    \ldots). 
  \end{equation}
  From lemma \ref{lem:etahom}, and the fact that the metric $h_0(r)$
  is $\epsilon^{-1}h_\epsilon$ with $\epsilon=\exp^{-r}$, the same
  holds true for $\eta(h(r)) - \eta_0(M,\theta)$.
  
  Moreover, the boundary contribution $B(g_{KE}, M_r)$ arising in
  definition \ref{defn:nu} of $\nu$ is the integral of a secondary
  class built from the curvature of $g_{KE}$ and has therefore a
  development of the same type as \eqref{eq:2:10}. The expression
  $$
  \nu(r) + 3\eta_0(M, \theta) = B(g_{KE}, M_r) - 3 ( \eta(h(r)) +
  \eta_0(M, \theta) )
  $$
  has then a development of the same kind.  Note that this
  expression is void of terms in $\exp^{kr}$ for $k > 0$ since we
  already know from definition \ref{defn:nu} and \cite{ob-mh1} that it
  converges when $r$ goes to infinity. As a result, the local boundary
  contribution necessarily cancels {\sl all} divergent terms, and adds
  (still local) convergent terms.  Identifying the constant terms we
  get eventually:
  \begin{displaymath}
    \nu (M) + 3 \eta_0(M, \theta) = \int_M P_\theta(R, \tau, \nabla R, \nabla
    \tau, \ldots) \,\theta \land d\theta  
  \end{displaymath}
  where $P_{\theta}$ is some pseudohermitian local tensorial
  invariant.  The invariance under the rescaling $\theta \to \lambda^2
  \theta$ shows that the polynomial $P_\theta$ must satisfy
  $$
  P_{\lambda^2 \theta} = \lambda^{-4 } P_\theta.
  $$
  The list of all possible expressions is easily established.
  Indeed, elementary invariant theory yields that such
  $U(1)$-invariant polynomials have to be sums of full contractions.
  Curvature $R$ and torsion $\tau$ (here we see the torsion $\tau$ as
  a tensor of type $\tau = A_{11} \theta^1 \otimes \theta^1$ using
  some coframe $\theta^1 $ of $T^{1,0} H$) are homogeneous of weight
  $-2$ with respect to the previous rescaling, while a covariant
  differentiation along $T$ decreases the weight by $2$, and an
  horizontal one by $1$ .  Following proposition 5.13 in
  \cite{Stanton}, we find that $P_\theta$ is a combination of
  \begin{equation}
    \label{eq:2:11} 
      \begin{gathered}
        R^2\ ,\ |\tau|^2 = |A_{11}|^2\ ,\ R_{,0} = dR(T)\ ,\ \Delta_H
        R\ , \\ \nabla^2_{0,1} \tau = A_{11, \bar 1 \bar 1} \ , \ 
        \nabla^2_{1,0} \bar\tau = A_{\bar 1 \bar 1, 11}.
    \end{gathered}
  \end{equation}
  Full divergences do not contribute after integration over $M$, so
  that one may forget the last four expressions, and the proof of
  lemma \ref{lem:2:3} is over.
\end{proof}

\subsection*{Computation of the constants}
We are left with the determination of $C_1$ and $C_2$ in lemma
\ref{lem:2:3}. This shall come from an explicit study of
left-invariant CR structures on the three sphere.

Choose a basis $(\alpha_1,\alpha_2,\alpha_3)$ of left-invariant 1-forms on $\cerc^3$,
such that $d\alpha_1=\alpha_2\land \alpha_3$, etc.  The $\eta$-invariant of the
left-invariant metric $\lambda_1^2\alpha_1^2+\lambda_2^2\alpha_2^2+\lambda_3^2\alpha_3^2$ has been
computed by Hitchin \cite[formula
(10)]{hitchin-bolletino}\footnote{There is a slight mistake in
  \cite{hitchin-bolletino} by a factor $2$, as can be seen by
  comparing the results in \cite{hitchin-bolletino} for the standard
  sphere to those of theorem \ref{th:ouy} below: one must find
  $\eta_0(\cerc^3, std)=\frac{2}{3}$ in the equation \eqref{eq:eta0li}
  below, rather than $\frac{4}{3}$ computed by \cite{hitchin-bolletino}.}:
\begin{equation}
  \label{eq:etali}
  \eta(\lambda_1^2\alpha_1^2+\lambda_2^2\alpha_2^2+\lambda_3^2\alpha_3^2) 
=\frac{2}{3} \left( \frac{s_1^3-4s_1s_2}{s_3}+ 9 \right)
\end{equation}
where the $s_i$ are the symmetric polynomials in the $\lambda_i^2$. As a result, we get
\begin{multline*}
  \eta(\alpha_1^2+\lambda_2^2\alpha_2^2+\lambda_3^2\alpha_3^2) \\ =
  \frac{2}{3\lambda_2^2\lambda_3^2} \big( \lambda_3^6 -
  (1+\lambda_2^2) \lambda_3^4 - (\lambda_2^4-3\lambda_2^2+1)
  \lambda_3^2 + (\lambda_2^6 - \lambda_2^4 -\lambda_2^2 + 1) \big)
\end{multline*}
and taking the constant term in the diabatic limit $\lambda_3 \to
\infty$ (\emph{i.e.} taking $\theta = \alpha_3$) leads to
\begin{equation}
  \label{eq:eta0li}
  \eta_0(\alpha_1^2+\lambda^2\alpha_2^2) 
    = \frac{2}{3\lambda^2} (-\lambda^4 + 3\lambda^2 -1) .
\end{equation}
On the other hand, the $\nu$-invariant can be estimated from the
$\mu$-invariant introduced by Burns and Epstein for embeddable CR
structures, or more generally CR manifolds with trivial holomorphic
part of the contact bundle \cite{BE88}: for the contact form
$\theta=\alpha_3$ and a metric $\gamma =
\lambda^{-1}(\alpha_1)^2+\lambda(\alpha_2)^2$, $\mu$ is calculated in
\cite[4.1.A]{BE88}. Since
\begin{equation}\label{eq:Rli}
  R = \frac{1+\lambda^2}{2\lambda}, \quad |\tau| = \frac{1-\lambda^2}{2\lambda}, 
\end{equation}
one has
\begin{equation*}
  \mu(\lambda^{-1}\alpha_1^2+\lambda\alpha_2^2) = - \frac{1}{16\pi^2}
  \int_{S^3} (4|\tau|^2 - R^2 ) \theta\land d\theta \notag = - 1 +
  \frac{3(1-\lambda^2)^2}{4\lambda^2} .
\end{equation*}
It is proved in \cite{ob-mh1} that, for a deformation of the standard
CR $3$-sphere, one has $\nu = 3\mu + 2$, and therefore
\begin{equation}
  \label{eq:nuli}
  \nu(\lambda^{-1}\alpha_1^2+\lambda\alpha_2^2) \ = \ -1 \ + 
\ \frac{9(1-\lambda^2)^2}{4\lambda^2} \ .
\end{equation}
From equations \eqref{eq:eta0li}, \eqref{eq:nuli} and \eqref{eq:Rli}
we deduce
\begin{equation*}
  (\nu+3\eta_0)(\lambda^{-1}\alpha_1^2+\lambda\alpha_2^2) \
  = \ \frac{(1+\lambda^2)^2}{4\lambda^2} \
  = \ \frac{1}{16\pi^2} \,\int_{S^3} R^2 \, \theta\land d\theta \ .
\end{equation*}
This yields $16\pi^2C_1=1$ and $C_2=0$ and the proof of theorem
\ref{eta0} is done. \qed

\smallskip

\begin{rem}\label{rem:ndnu}
  From Theorem \ref{eta0}, we see that $\displaystyle
  -3\eta_0+\frac{1}{16\pi^2}\int_MR^2\theta\land d\theta$ is a CR
  invariant. This fact can be proved directly: standard calculations
  in pseudohermitian geometry lead easily to the conclusion that it is
  invariant under conformal transformations $\theta\to f\theta$.
  
  This provides an alternative (and independent) definition of the
  $\nu$-invariant.  The latest is clearly simpler than the one
  explained in section \ref{sec:nu}: this is useful for computations
  and theoretical aspects, in particular the relation with the
  $\eta$-invariant of the contact operator $D*$ on vertical 2-forms,
  as we shall see in the following sections. On the other hand, very
  important for the applications is the fact that $\nu$ arises as a
  boundary term in the integral of characteristic classes (see for
  example corollary \ref{cor:3}), and this can be obtained only
  through the first definition and the work done in \cite{ob-mh1}.
  
  One may also think that this remark could serve as a basis for
  defining a version of $ \nu$ in higher dimensions, by looking for
  local corrections of $\eta_0$ that would lead to a CR invariant.
  However, this seems a very difficult task, as the range of possible
  terms of the right weight is in general much larger than in
  \eqref{eq:2:11}, even in the next relevant dimension $7$.
\end{rem}

\bigskip

\section{Computation of the invariant on Seifert manifolds}\label{sec:Sfrt}

This section is devoted to explicit computations of the
$\nu$-invariant on $\cerc^1$-invariant CR manifolds of dimension $3$.
Although certainly a digression from our main route towards Theorems
\ref{th:D*eta0} and \ref{maintheo2}, this appears as a nice direct
application of the results obtained in the previous section. We have
thus chosen to interrupt the pace of our proofs, and to offer this
section as a refreshing intermezzo before the analytical
technicalities that will follow.

\smallskip

We first describe our family of spherical $3$-dimensional compact
strictly pseudoconvex CR manifolds in greater detail.

\begin{defi}\label{defn:3:1}
  A CR-Seifert manifold is a $3$-dimensional compact manifold endowed
  with both a pseudoconvex CR structure $(H,J)$ and a Seifert
  structure, that are compatible in the following sense: the circle
  action $\varphi : \cerc^1 \to \mathrm{Diff}(M)$ preserves the CR
  structure and is generated by a Reeb field $T$.
\end{defi}

Any $\cerc^1$-invariant CR structure admits a $\cerc^1$-invariant
contact form $\theta$ if the manifold is orientable (this is proved in
\cite{ka-tsu}). Moreover it is easily proved that that existence of a
Reeb field $T$ (defined by $\theta(T)=1$ and $\iota_Td\theta = 0$)
satisfying $\varphi_*(\frac{d}{dt}) = T$ and $\mathcal{L}_T\theta =
0$, $\mathcal{L}_T J= 0$, is equivalent to the existence of a
locally free action of $\cerc^1$ whose (never vanishing) infinitesimal
generator preserves $H$ and $J$ and is transverse everywhere to $H$.
Hence, our CR-Seifert manifolds could also be called \emph{transverse
  $\cerc^1$-invariant CR manifolds}; note moreover that there exists a
much larger class of $\cerc^1$-invariant CR manifolds, with the
infinitesimal generator being sometimes tangent to the contact
distribution \cite{ka-tsu,lutz}.

As we do not assume the action to be free but only locally free, the
quotient space $\Sigma=M/\cerc^1$ is a surface with possibly conical
singularities.  Each CR-Seifert manifold is then an orbifold bundle
over the compact Riemannian orbifold surface $\Sigma$. If $\Sigma$ is
such a surface, endowed with a complex structure, orbifold
$\cerc^1$-bundles are classified by their (rational) degrees $d$.
Singularities of the bundle are located above the singularities of
$\Sigma$ in such a way that the resulting $3$-manifold is smooth: if
the local fundamental group is $\setZ/\alpha\setZ$ ($\alpha\in\NM^*$),
a generator acts on a local chart around $p$ of the basis manifold as
$\exp^{i\frac{2\rho\pi}{\alpha}}$ and on the fiber as
$\exp^{i\frac{2\pi\beta}{\alpha}}$ with $\rho$ and $\beta$ prime to
$\alpha$ (the extra parameter $\rho$ may seem pointless as it is 
always possible to reduce oneself to two parameters by taking $\rho'=1$ 
and $\beta'=\beta\rho^{-1}$ mod. $\alpha$, but this extended description
will prove useful when specializing our computations to the case of lens 
spaces in section \ref{sec:coro}) . 
Any choice of equivariant connection $1$-form $\theta$ on
$M$ endows it with an invariant CR structure, $H$ being chosen as the
horizontal space for the connection and $J$ being pulled back from
the base. It is strictly pseudoconvex if $d<0$. The interested reader
is referred to \cite{Nicolaescu2} for a very readable account on
orbifold bundles over orbifold surfaces. Note moreover that one has
\[ \int_M \theta\land d\theta = -4 \pi^2 d, \]
and that the metric $\gamma = d\theta(\cdot,J\cdot)$ projects
downwards to a metric on $\Sigma$ of volume
\[ \int_{\Sigma} d\theta = -2\pi d,\]
(see \cite{Nicolaescu2} again for integration of forms over
orbifolds). Its curvature $R$ equals the Tanaka-Webster curvature of
$(M,\theta)$ and Gauss-Bonnet reads
\[ \int_\Sigma R\, d\theta = 2\pi\,\chi , \]
where $\chi$ is the (rational) Euler characteristic of $\Sigma$.

\medskip

\subsection*{Computations in constant curvature} 
In the first half of this section, we moreover assume that $\gamma$
has constant curvature $R$. In this case, the CR structure is
spherical, that is $M$ is locally isomorphic to the standard
$3$-sphere. Conversely, it is known that spherical CR-Seifert
manifolds are exactly those of constant Tanaka-Webster curvature $R$,
except if the base is a sphere, see for instance \cite{Belgun}.

The computations now rely on the explicit derivation of the
$\eta$-invariant of (orbifold) circle bundles over (orbifold)
Riemannian surfaces with constant curvature that have been done by
Komuro \cite{komuro} and more generally by Ouyang \cite{ouyang}.  In
our conventions and notations, their results read:

\begin{theo}[Ouyang]\label{th:ouy}
  The $\eta$-invariant of the metric $t^2\,\theta^2 + \gamma$ on
  $M$ is equal to
  \begin{displaymath}
    \frac{1}{3} \left( d + 3 + 2d \left( \frac{\pi t^2}{V}\chi - 
        \frac{\pi^2 t^4}{V^2}d^2 \right)\right) + 4 \sum_{j=1}^{p}
    s(\alpha_j,\rho_j, \gamma_j ) ,
  \end{displaymath}
  where $s(\alpha,\rho,\gamma)=\frac{1}{4\alpha}\sum_{k=1}^{\alpha
    -1} \cot (\frac{k\rho\pi}{\alpha})\cot (\frac{k\beta
    \pi}{\alpha})$ is the classical Rademacher-Dedekind sum.
\end{theo}

We can now proceed to the computation of $\nu$ in the constant
curvature case.  We have to show Corollary \ref{maintheo3}, which we restate
here:

\begin{cor}\label{nuSfrt}
  Let $M$ be a compact $\cerc^1$-orbifold bundle of rational degree
  $d<0$ over a compact orbifold surface $\Sigma$ of constant curvature
  and rational Euler characteristic $\chi$. Then,
\begin{equation}\label{eq:Seifert-calcul1}
\nu(M) = -d - 3 - \frac{\chi^2}{4d} - 12 \sum_{j=1}^{p}
s(\alpha_j,\rho_j, \beta_j). 
\end{equation}
\end{cor}

Let us remark that the $\nu$-invariant depends only on the topology
for this class of CR manifolds, and not, for instance, on the complex
structure of $\Sigma$.  This is {\it a priori} known, since the
gradient of $\nu$ is the Cartan curvature \cite[Theorem 8.1]{ob-mh1},
which vanishes for spherical CR manifolds.

\begin{proof}
  According to Theorem \ref{eta0}, the $\nu$-invariant is given by
  adding a local term to the renormalized $\eta$-invariant.  On
  $\cerc^1$-invariant CR manifolds with constant curvature, the
  renormalized invariant is easily read from Ouyang's Theorem
  \ref{th:ouy} above:
  \begin{equation}
    \label{eq:2}
    \eta_0 \ = \ 1 + \frac{d}{3} + 4\,\sum_{j=1}^{p} 
  s(\alpha_j,\rho_j, \beta_j) \ .
  \end{equation}
  Moreover, the integral term is just
  \[ \frac{1}{16\pi^2}  \,\int_{M} R^2 \theta\land d\theta \ = 
  \ \frac{-4\pi^2d\left( - \frac{\chi}{d}\right)^2}{16\pi^2} \ = \ -
  \,\frac{\chi^2}{4d} \,,
  \]
  which shows also Theorem \ref{maintheo3} in the constant curvature
  case.
\end{proof}

\medskip

\begin{rem}
  Corollary \ref{maintheo3} can also be obtained by direct calculation
  from the original definition of $\nu$ and Ouyang's formula. 
  Indeed the asymptotically
  Kähler-Einstein metric $g_{KE}$ on $[r_0,+\infty[× M$ can be handled
  with bare hands in this simple situation, and the boundary
  contribution counterbalancing the divergence of the sequence of
  $\eta$-invariants can be explicitly derived. Putting together
  Ouyang's theorem \ref{th:ouy} and these local computations yield the
  value of $\nu$, see \cite{mh-vol} for similar computations. This is
  of course a painful method, but it is still a reasonably simple case
  where the cancellation of divergences by local terms can be observed
  in detail.
\end{rem}

\medskip

\subsection*{Extension to cases of non-constant curvature}
We now extend the computations of $\nu$ to an (almost) complete proof
of theorem \ref{maintheo1}. It is shown in \cite{ka-tsu,lutz} that
there always exist a unique (up to equivalence) transverse
$\cerc^1$-contact form on an orientable Seifert manifold (careful:
this might be wrong for a non-transverse action). Given the natural
contact form that fixes the length of the regular fibers to $2\pi$,
the choice of a CR structure is then equivalent to the choice of a
downwards orbifold Riemannian metric $\gamma$ of fixed volume
$d\theta$, and this metric might or might not be of constant
curvature.

In case the base is smooth (no orbifold singularities), it is known
that the adiabatic limit $\eta_{\mathrm{ad}}$ does not depend on the
underlying metric on $\Sigma$, see \emph{e.g.} \cite{Zhang}. As one
can always find a constant curvature metric of volume $d\theta$ (easy
consequence of Moser's lemma on volume forms), the previous formula
\eqref{eq:2} for $\eta_0= \eta_{ad}$ applies.
Then Theorem \ref{eta0} enables to conclude that
\begin{equation}\label{eq:Sfrt-general}
\nu(M) = -d - 3  - 12 \sum_{j=1}^{p} s(\alpha_j,\rho_j, \beta_j)
+  \frac{1}{8\pi}\, \int_\Sigma R^2\, d\theta . 
\end{equation}
If orbifolds singularities are present, it is known that every
orbifold surface has a constant curvature metric, except some
exceptional cases on the sphere described in \cite{Belgun2}.  As the
set of compatible complex structures with a given contact structure is
contractible, this means that, except on the exceptional cases we have
just alluded to, it suffices to check the following:

\begin{lem}\label{lem:vareta0}
The variations of $\eta_0$ with respect to the complex structure vanish 
when the torsion is zero.
\end{lem} 

\begin{proof}
  From Theorem \ref{eta0}, $\eta_0$ has the same variation as
  $$
  -\frac{\nu}{3} + \frac{1}{48\pi^2}\,\int_M R^2 \,\theta \land
  d\theta .
  $$
  The variation of $\nu$ with respect to $J$ has been computed in
  \cite[Theorem 8.1]{ob-mh1}, namely
  \begin{equation}
    \label{eq:3}
    \frac{d \nu}{dJ} = \frac{-3}{8\pi^2} \int_M \langle Q_J, \dot
    J\rangle \,
    \theta \land d \theta\, , 
  \end{equation}
  where $ Q_J = i {Q_1}^{\bar 1} \theta^1 \otimes Z_{\bar 1}
  - i {Q_{\bar 1}}^1 \theta^{\bar 1} \otimes Z_1 \in \End
  (H)$ is Cartan's tensor. Its expression in term of derivatives of
  Tanaka-Webster curvature and torsion is given by
  \begin{equation}
    \label{eq:4}
    {Q_1}^{\bar 1} = \frac{1}{6} {R_{,1}}^{\bar 1} +
    \frac{i}{2} R {A_1}^{\bar 1} - {A_1}^{\bar 1}{}_{,\,0} -
    \frac{2i}{3} {A_1}^{\bar 1}{}_{,\bar 1}{}^{\bar 1} \,.
  \end{equation}
  On the other hand the variation of the Tanaka-Webster curvature is
  computed e.g in \cite[(2.20)]{cheng-lee}, and is given by
  \begin{equation}
    \label{eq:5}
    \dot R = i(E_1^{\ \bar 1}{}_{,\bar 1}{}^1 - E_{\bar
    1}^{\ 1}{}_{,1}{}^{\bar 1}) - ({A_1}^{\bar 1} E_{\bar
    1}^{\ 1} + A_{\bar 1}^{\ 1} E_1^{\ \bar 1})\,, 
  \end{equation}
  where 
  \begin{equation}
    \label{eq:6}
    \dot J= 2 E_1^{\ \bar 1} \theta^1 \otimes Z_{ \bar 1} + 2
    E_{\bar 1}^{\ 1} \theta^{\bar 1} \otimes Z_{\bar 1}. 
  \end{equation}
  Putting everything together and integrating by parts shows that, in
  vanishing torsion, $\eta_0$ does not depend on the complex structure
  as needed.
\end{proof}

\begin{rem}
  This computations of variations may be seen as an alternative mean to
  determine the constant $C_1= \frac{1}{16\pi^2}$ in Lemma
  \ref{lem:2:3}, independently of the computations of examples done in
  section \ref{sec:2:3}. Moreover it shows that $\eta_0$ is
  independent of $J$ whenever the torsion vanishes, \emph{without any
  assumption on the quotient structure of $M$ by the Reeb flow}.
  This last fact will be used in section 9.
\end{rem}

In the remaining exceptional cases over $\cerc^2$ described in
\cite{Belgun2}, the results stay the same but the proof above 
does not apply anymore and one has to rely on
a different technique: this will be done below in section
\ref{sec:decompostion-sp-d}. 

\bigskip

\section{The contact complex and the diabatic limit.}
\label{sec:4:1}

Theorem \ref{eta0} gives a simple formula relating the $\nu$-invariant
and the renormalized $\eta$-invariant $\eta_0$ of the contact-rescaling.
According to \eqref{eq:1}, $\eta_0$ coincides with the adiabatic limit of
$\eta$ in the case the CR manifold has vanishing torsion, and this
enables computations, for explicit expressions of the adiabatic limit
are known in a number of cases. But a deeper question is to relate
directly the $\nu$-invariant to the geometry and spectral theory of the
CR or pseudohermitian manifold.

In the sequel we shall consider a natural $\eta$-invariant arising in
pseudohermitian geometry. One actually knows by \cite{Rumin00} a
candidate for this, coming from the contact-de Rham complex. We shall
briefly recall its construction in dimension $3$ and its relation with
the diabatic limit.

\medskip

Let $M$ be a $3$-dimensional contact manifold and $H$ its contact
distribution. We denote by $\Omega^* H$ the space of horizontal forms,
\emph{i.e.} the space of sections of the alternating algebra over the
dual of the bundle $H$. Let also $\Omega^* V$ be the subspace of
vertical forms on $M$, by which we mean ``true'' forms in $\Omega^*M$
vanishing on $H$. Equivalently, one has $\Omega^* V= \{ \theta \land
\alpha\} = \theta \land \Omega^* H$ for any local choice of contact
form $\theta$.  The contact-de Rham complex is then the following:
\begin{equation}
  \label{eq:7}
  C^{\infty}(M) \stackrel{d_H}{\longrightarrow} \Omega^1H
  \stackrel{D}{\longrightarrow} \Omega^2 V 
   \stackrel{d_H}{\longrightarrow} \Omega^3 M,
\end{equation}
where for $f \in C^\infty (M)$, $d_H f \in \Omega^1H$ stands for the
restriction of $df$ to $H$, while
$$
d_H : \Omega^2 V \to \Omega^3 M
$$
is just de Rham's differential restricted to $\Omega^2 V$ in
$\Omega^2 M$, and $D$ is defined as follows: since $d$ induces an
isomorphism
$$
d_0 : \Omega^ 1 V \to \Omega^ 2 H \ \text{ with } \ d_0(f
\theta ) = f d \theta_{\mid \Lambda^2H},
$$
then any $\alpha$ in $\Omega^1 H$ admits a unique extension $\ell
(\alpha)$ in $\Omega^1M $ such that $d \ell( \alpha)$ belongs to
$\Omega^2 V$; namely, given any initial extension $\overline \alpha$
of $\alpha$, one has
\begin{equation}
  \label{eq:8}
  \ell (\alpha) = \overline \alpha - d_0^{-1} (d \overline
  \alpha)_{\mid \Lambda^2H}. 
\end{equation}
We then define
\begin{equation}
  \label{eq:9}
  D\alpha= d \ell (\alpha). 
\end{equation}
This differential $D$ is a second order operator, since the lifting
$\ell : \Omega^ 1 H \to \Omega^ 1 M$ is a first order one.
Moreover one sees easily that $\ell$ induces an homotopy equivalence
between the contact and de Rham complexes, together with the natural
restrictions, and the retraction $\ell' : \Omega^2 M \longrightarrow
\Omega^2 V$ defined by
\begin{displaymath}
  \ell' (\alpha)  = \alpha - d d_0^{-1} \alpha_{\mid \Lambda^2 H}.
\end{displaymath}
From now on we will suppose moreover that the contact manifold $M$ is
endowed with a strictly pseudoconvex CR structure $J$, together with
some choice of contact form $\theta$. We consider the
contact-rescaling sequence of metrics of \eqref{eq:2:1}
\begin{displaymath}
  h_0(r)  = \exp^{2r} \theta^2 + \exp^r d \theta (\cdot, J\cdot).
\end{displaymath}
Let $\varepsilon = \exp^{-r}$, as before, and define
\begin{equation}
  \label{eq:10}
  g_\varepsilon = \varepsilon^{-2} \theta^2 + \varepsilon^{-1} d
  \theta (\cdot, J 
  \cdot) = h_0(r).
\end{equation}
This metric induces an orthogonal splitting $ TM = H \oplus \setR T$
where $T$ is the Reeb field of $\theta$, and one can identifies
$\Omega^1 H$ with ``true'' $1$-forms on $M$ vanishing on $T$.
Observing that Hodge $*$-operator exchanges $\Omega^1 H$ and $\Omega^2
V$ and one can consider $D*$ acting on closed vertical $2$-forms
$\Omega^2_D V = \Omega^2 V \cap \im D$.

Following \cite[Theorem 4.14]{aps1}, we define the boundary operator for
the signature attached to the Riemannian metric $g_\epsilon$ as
$$
S_\varepsilon = (-1)^p (*_\varepsilon d - d*_\varepsilon),
$$
acting on $\Omega^{2p} M = C^\infty M \oplus \Omega^2 M $. As
observed in \cite[Prop 4.20]{aps1}, one may remove some spectral
symmetry, and its $\eta$-function
\begin{equation}
  \label{eq:11}
  \eta(S_\varepsilon)(s) = \tr^* (S_\varepsilon
  |S_\varepsilon|^{-(s+1)}) = \sum_{\lambda_i \in \spec(S_\varepsilon)
  \setminus \{0\}
  } \frac{\lambda_i}{|\lambda_i |^{s+1}}
\end{equation}
actually coincides with that of $d *_\varepsilon$ when restricted to
$\Omega^2_d M = \Omega^2 M \cap \im d$. Note that we have used $\tr^*$
to denote a trace taken outside the $0$-eigenspace. In the same vein,
the notation $\spec^*$ used below will denote a spectrum \emph{where
  the $0$-eigenvalue has been removed}.

From \cite[p. 74]{aps3} or \cite[Chap. 1.10]{Gilkey}, the series
\eqref{eq:11} is absolutely convergent for $\Re s > 3$ and has a
meromorphic extension to $\C$, with possibly simple poles at $s = 3 -
n$, $n \in \N$. By Atiyah-Patodi-Singer's theorem \cite{aps1},
$\eta(S_\varepsilon)(s)$ is actually regular at $s=0$ and its value
there is called the $\eta$-invariant of $(M, g_\varepsilon)$.
Similarly, an $\eta$-function and its value at $0$ can be defined for
the operator $D*$ in dimension $3$. This mainly follows by applying
the same ideas, but with the adequate symbolic calculus for
hypoelliptic operators, see section \ref{sec:eta+}. 

\smallskip

In order to compare them, let us now compute $d *_\varepsilon$ and
$D*_\varepsilon$ using the decomposition of $\Omega^2 M$ into vertical
and horizontal $2$-forms:
\begin{displaymath}
  \alpha = \theta \land \alpha_T + \alpha_H,
\end{displaymath}
with $\alpha_T \in \Omega^1 H$, $\alpha_H \in \Omega^2 H$.  From
\eqref{eq:10} one sees that
$$
*_\varepsilon \alpha= \varepsilon *_H \alpha_T + \theta \land *_H
\alpha_H
$$
where $*_H$ denotes the induced Hodge duality on $H$. In matrix
form, one gets
\begin{equation}
  \label{eq:12}
  d *_\varepsilon = 
  \begin{pmatrix}
    \varepsilon\,\mathcal{L}_T *_H & -d_H *_H\\
    \varepsilon\, d_H *_H & 1
  \end{pmatrix},
\end{equation}
where $\mathcal{L}_T$ is the Lie derivative along $T$.

Using \eqref{eq:8} and \eqref{eq:9} one finds that $ \ell(\beta) =
\beta - (*_H d_H \beta ) \theta $ on $ \Omega^1 H$, so that $D \beta =
\theta \land (\mathcal{L}_T + d_H *_H d_H ) \beta$, and hence
\begin{equation}
  \label{eq:13}
    D *_\varepsilon (\theta \land \alpha_T) = \varepsilon \theta\land
    (\mathcal{L}_T + d_H *_H d_H )*_H \alpha_T
\end{equation}
on $\Omega^2 V = \theta \land \Omega^1 H$.

The whole spectrum of $D *_\varepsilon = \varepsilon D*_1$ then
collapses at speed $\varepsilon$ in the diabatic limit $\varepsilon
\to 0$, whereas part of the spectrum of $d *_\varepsilon$ is
not collapsing: for instance $(d *_\varepsilon) (d\theta )= d \theta$.
Hence the diabatic behaviour of the \emph{whole} spectrum of $d
*_\varepsilon$ cannot be related to $D*_\varepsilon$ alone, and indeed
only the collapsing spectra are related.  This shows up in the
following formulas, which are direct consequences of \eqref{eq:12} and
\eqref{eq:13}, or even more directly from the definitions \eqref{eq:8}
and \eqref{eq:9} of $\ell$ and $D$. If $P_\varepsilon=
\varepsilon^{-1}d *_\varepsilon$,
\begin{equation}
  \label{eq:14}
  \begin{split}
    P_\varepsilon = \varepsilon^{-1} d *_\varepsilon & =
  \begin{pmatrix}
    D*_1 & 0\\
    0 & 0
  \end{pmatrix} + 
  \begin{pmatrix}
    -(d_H *_H)^2 & - \varepsilon^{-1}d_H *_H\\
    d_H *_H & \varepsilon^{-1}
  \end{pmatrix} \\
  & = \Pi_{\Omega^2 V} (D *_1) \Pi_{\Omega^2 V} + \varepsilon
  P_\varepsilon \Pi_{\Omega^2 H} P_\varepsilon.
  \end{split}
\end{equation}
It follows that in the diabatic limit $\varepsilon \to 0$ all
the eventually bounded spectrum of $P_\varepsilon= \varepsilon^{-1}d
*_\varepsilon$ converges, at least weakly, towards the spectrum of
$D*_1$.  Actually its turns out that this spectral convergence is
uniform over bounded intervals, as a consequence of the uniform
convergence in the diabatic limit of the resolvents $(\lambda -
P_\varepsilon)^{-1}$ on $\ker d$ towards $(\lambda - D*_1)^{-1}$, for
$\lambda \in \C \setminus \R$ \cite[theorem 3.6]{Rumin00}.

Such a spectral convergence is unfortunately only a first step in the
study of a global spectral invariant like $\eta$.  To illustrate this,
recall that by \cite{Bi-Fr} an equivalent expression of the Riemannian
$\eta$-invariant is given by
\begin{equation}
  \label{eq:15}
  \eta(P_\varepsilon)(0) = \pi^{-1/2}\int_0^\infty \tr
  \bigr(P_\varepsilon \exp^{-t P_\varepsilon^2} \bigl) \frac{dt}{\sqrt t}.
\end{equation}
Now by \cite[Theorem 7.1]{Rumin00} the following global trace convergence
holds
$$
\tr (P_\varepsilon \exp^{-t P_\varepsilon^2} ) \longrightarrow \tr
( D* \exp^{-t DD^*}),
$$
when $\varepsilon $ goes to $0$, but uniformly on $t$ \emph{only
  for $t\geq t_0 > 0$}. It cannot be true for small $t$ since the
$\eta$-invariants and the integrals \eqref{eq:15} diverge in the
diabatic limit (although one knows by transgression formulas that
these divergences of $\eta(P_\varepsilon)(0)$ are given by local
expressions).  From the analytic viewpoint, these divergences are
rooted in the transition from elliptic towards hypoelliptic operators,
that cannot be uniform in all $(t, \varepsilon)$ regimes. For
instance, the asymptotic spectral densities (Weyl's laws), or the
powers of $t$ occurring in the asymptotic expansions of the heat
kernels for $t\to 0$ are not the same for the elliptic $P_\varepsilon$
and the hypoelliptic $D*$.  However it is possible, as is usual in
such asymptotic spectral problems, that the divergences occurring in
the $(d*_\varepsilon, D*)$ transition when $\varepsilon$ and $t$ go to
$0$, are ruled again by local expressions in the curvature, see also
Remark \ref{rem:conjecture}.  This would provide directly a relation
like \eqref{eq:etanugeneral} between the finite part $\eta_0$ of
$\eta(P_\varepsilon)$ in the diabatic limit and the contact
$\eta$-invariant $\eta(D*)$.  Unfortunately, the techniques used in
\cite{Rumin00} cannot handle these problems in the general case. The
analysis can however be done in the particular case of CR-Seifert
manifolds, and we will now restrict ourselves to this case.

\bigskip

\section{Spectral analysis on Seifert manifolds.}
\label{sec:metrics-h_0r-from}

As explained above, we will now deal with CR-manifolds endowed with
both a Seifert and a CR structure compatible in the sense that the
circle action $\varphi : \cerc^1 \to \mathrm{Diff}(M)$
preserves the CR structure $(H, J)$ and is generated by a Reeb field
$T$. An invariant contact form $\theta$ has then been chosen, and we
note that in this section, opposite to section \ref{sec:Sfrt}, we will
never assume the Webster curvature to be constant.

\medskip

The circle action allows to perform a Fourier decomposition of
functions or forms, inside $M$ and without referring to the quotient
structure. For instance, given $n \in \Z$ and $f \in C^0(M)$, its
$n$-th component is the function on $M$ defined by
$$
\pi_n f = \frac{1}{2\pi}\int_0^{2\pi} \exp^{-int} (f\circ
\varphi_t) dt.
$$
It satisfies $(\pi_n f) \circ \varphi_t = \exp^{int} (\pi_n f)$, so
that $\mathcal{L}_T (\pi_n f) = in \pi_n f$ on $C^1(M)$. The
projections $\pi_n$ preserve and are clearly bounded on all $C^p(M)$,
$L^p(M)$ or Sobolev spaces. Moreover, the Hilbert sum of all $\pi_n$
for $n$ in $\Z$ is the identity on $L^2(M)$.  Last, this circle action
preserves all structures and operators related to the above choice of
contact form, so that we will be able to split their spectra into
Fourier components.

\smallskip

We can now study the spectral aspects of the contact rescaling
$g_\varepsilon$ in \eqref{eq:10} on a CR-Seifert manifold $M$.  Of
course the adiabatic limit exists in this situation, and has already
been much studied, see \emph{e.g.} \cite{bismut-cheeger, dai}, but we
will need a different approach here, focusing on the diabatic
behaviour of $d*_\varepsilon$ and $\eta(d*_\varepsilon)$, as related
to the spectrum of $D*$ and its $\eta$-invariant.

One computes easily the Laplacian on $\Omega^2 M$, relatively to the
splitting
$$\Omega^2 M = \theta \land \Omega^1H \oplus \Omega^2H, $$
namely
\begin{equation}
  \label{eq:16}
  \Delta_\varepsilon = 
  \begin{pmatrix}
    \varepsilon \Delta_H - \varepsilon^2 T^2 & - d_H *_H\\
    \varepsilon d_H *_H & 1 + \varepsilon \Delta_H - \varepsilon^2 T^2
  \end{pmatrix},
\end{equation}
where $\Delta_H = d_H \delta_H + \delta_H d_H$ is the horizontal
Laplacian (not to be confused with the contact Laplacian introduced in
\cite{Rumin94,Rumin00}), $T$ denotes here the Lie derivative along
$T$, and we have used that $T^* = -T$ and $[T, \delta_H] =0$ since $T$
is a Killing Reeb field on the CR-Seifert manifold. We observe from
\eqref{eq:12} that the non diagonal part of $\Delta_\varepsilon$ is
the same as that of $d *_\varepsilon$, so that
$$
\Delta_\varepsilon = d *_\varepsilon \ + \ \varepsilon
\begin{pmatrix}
  \Delta_H + T *_H & 0 \\
  0 & \Delta_H
\end{pmatrix} \ - \ \varepsilon^2 T^2.
$$
When studying spectral asymmetry, we restrict ourselves to the
subspace $\Omega^2_d M =\im d$ of $\Omega^2M$, on which
$\Delta_\varepsilon = (d*_\varepsilon)^2$. We get therefore the
following expression relating pairwise commuting operators:
\begin{equation}
  \label{eq:17}
   (d *_\varepsilon)^2 = (d*_\varepsilon) +  (\varepsilon K) 
  - \varepsilon^2 T^2 \,, 
\end{equation}
with
$$
K = \begin{pmatrix}
  \Delta_H + T *_H & 0 \\
  0 & \Delta_H
\end{pmatrix} .$$ 
Therefore if $\alpha \in \Omega^2_dM \setminus \{0\}$ satisfies
\begin{equation}
  \label{eq:18}
  (d *_\varepsilon )\alpha = \lambda_\varepsilon \alpha \ , \quad
  K\alpha = k\alpha \quad \text{and} \quad T^2\alpha  = -n^2\alpha ,
\end{equation}
for $\lambda_\varepsilon$ a non-zero eigenvalue of $d*_\epsilon$, then
\begin{equation}
  \label{eq:19}
    \lambda_\varepsilon + \varepsilon k 
  +\varepsilon^2 n^2 = \lambda_\varepsilon^2 \not= 0\ ,
\end{equation}
and, necessarily,
\begin{equation}
  \label{eq:20}
  \lambda_\varepsilon = \lambda_\varepsilon^+ \ \mathrm{or}\
  \lambda_\varepsilon^- 
  \quad \mathrm{with } \quad\lambda_\varepsilon^{±} = \frac{1
    ± \sqrt{1+ 4 \varepsilon (k + 4 \varepsilon n^2)}}{2}. 
\end{equation}
Hence the spectrum of $d *_\varepsilon$ splits in two families which
behave differently in the diabatic limit $\varepsilon \to 0$.
Eigenvalues of type $\lambda_\varepsilon^-$ all collapse, while those
of type $\lambda_\varepsilon^+ $ all converge to $1$. According to the
general results of \cite{Rumin00} discussed in section~\ref{sec:4:1},
only eigenvalues of type $\lambda^-_\varepsilon$ are related to $D*$,
after rescaling by $\varepsilon^{-1}$.

\smallskip

The previous eigenvalue equation \eqref{eq:20} is only a necessary
condition and we have to determine which of the possible
$\lambda_\varepsilon^{±}$ are effectively present in
$\spec(d*_\varepsilon)$ and to compute their multiplicities.  To do
this, we use the splitting the induced by the choice of the Reeb
field: suppose $\alpha = \theta \land \alpha_T + \alpha_H$ is a
$2$-form in the image of $d$. By \eqref{eq:12}, the system $(d
*_\varepsilon) \alpha = \lambda_\varepsilon \alpha$ is
\begin{align}
  \label{eq:21}
  (\lambda_\varepsilon - \varepsilon T*_H)\,\alpha_T & = -d_H *_H
  \alpha_H\\
  \label{eq:22}
  (\lambda_\varepsilon - 1)\,\alpha_H & = \varepsilon\, d_H *_H
  \alpha_T.
\end{align}
Suppose now that
\begin{equation}
  \label{eq:17bis}
  (d *_\varepsilon )\alpha = \lambda_\varepsilon \alpha \ , \quad
  K\alpha = k\alpha \quad \text{and} \quad T^2\alpha  = -n^2\alpha .
\end{equation}
Then we observe that $*_H = -J $ on $\Omega^1 H$ and $(T*_H)^2 =
-T^2 = n^2$.  Therefore \eqref{eq:21} gives
\begin{equation}
  \label{eq:23}
  (\lambda_\varepsilon^2 - \varepsilon^2 n^2) \alpha_T = -
  (\lambda_\varepsilon  + \varepsilon T*_H)\ d_H *_H \alpha_H,
\end{equation}
so that $\alpha_H$ determines uniquely $\alpha_T$ when
$\lambda_\varepsilon^2 \not= \varepsilon^2 n^2$.  A first (quite
large) part of the non-zero spectrum is then handled as follows.

\begin{prop}\label{prop:4:5}
  $\bullet$ Forms $\alpha = \theta \land \alpha_T + \alpha_H$ in
  $\Omega_d^2M$ satisfying
  \begin{equation}
    \label{eq:17bisagain}
    (d *_\varepsilon )\alpha = \lambda_\varepsilon^+ \alpha \ , \quad
    K\alpha = k\alpha \quad \text{and} \quad T^2\alpha  = -n^2\alpha 
  \end{equation}
  are in one-to-one linear correspondence with forms $\alpha_H$ in
  $\Omega^2H$ satisfying
  \begin{equation}
    \label{eq:alphaH}
    \Delta_H \alpha_H = k\alpha_H \quad \text{and} \quad T^2\alpha_H
    = -n^2\alpha_H . 
  \end{equation}
  $\bullet$ Forms $\alpha = \theta \land \alpha_T + \alpha_H$ in
  $\Omega_d^2M$ satisfying
  \begin{equation}
    \label{eq:17bisagainagain}
    (d *_\varepsilon )\alpha = \lambda_\varepsilon^- \alpha \ , \quad
    K\alpha = k\alpha \quad \text{and} \quad T^2\alpha  = -n^2\alpha 
  \end{equation}
  such that $(\lambda_\varepsilon^-)^2\not= \epsilon^2n^2$ are in
  one-to-one linear correspondence with forms $\alpha_H$ in
  $\Omega^2H$ satisfying
  \begin{equation}
    \label{eq:alphaHbis}
    \Delta_H \alpha_H = k\alpha_H \quad \text{and} \quad T^2\alpha_H
    = -n^2\alpha_H  
  \end{equation}
  with $k\not= |n|$.
\end{prop}

\begin{proof}
  Note first that, for any eigenvector $\alpha$ of $d*_\epsilon$
  satisfying either \eqref{eq:17bisagain} or
  \eqref{eq:17bisagainagain}, one may have $(\lambda_\varepsilon)^2 =
  \epsilon^2n^2$ only if \eqref{eq:17bisagainagain} holds. Hence, 
  in the positive case, one always has $\alpha_H\not= 0$, and, as a result,
  $\Delta_H \alpha_H = k \alpha_H$, $k$ is necessarily non-negative
  and $T^2\alpha_H = -n^2\alpha_H$. In the negative case, the same
  holds only if $(\lambda_\varepsilon^-)^2\not= \epsilon^2n^2$, and
  \eqref{eq:19} shows that this is equivalent to $k\not= |n|$.
  
  \smallskip
  
  Conversely, suppose now given $\alpha_H$, $n$, $k$,
  $\lambda_\varepsilon$ as needed. From \eqref{eq:23}, one defines
  $$
  \alpha_T = -(\lambda_\varepsilon^2 - \varepsilon^2 n^2)^{-1}
  (\lambda_\varepsilon + \varepsilon T*_H)\ d_H *_H \alpha_H,
  $$
  which satisfies \eqref{eq:21}. To check \eqref{eq:22}, recall
  that
  $$
  \delta_H = - *_H d_H *_H \quad \text{and} \quad d_H ^2 = - L T =
  -TL,
  $$
  (the last equation being a consequence of $d^2= 0$ see e.g.
  \cite[p. 415]{Rumin00} with $L(f)=f\,d\theta$).  One finds
\begin{align*}
  (\lambda_\varepsilon^2 - \varepsilon^2 n^2) d_H *_H \alpha_T & =
  (\lambda_\varepsilon d_H \delta_H \alpha_H - \varepsilon d_H^2
  T*_H \alpha_H)\\
  & = (\lambda_\varepsilon \Delta_H - \varepsilon T^2) \alpha_H \\
  & = (\lambda_\varepsilon k + \varepsilon n^2) \alpha_H.
\end{align*}
The eigenvalue equation \eqref{eq:19} then easily leads to
\eqref{eq:22}.
\end{proof}

For later use, note that the choice $(k,n)=(0,0)$ in the positive case
leads to $\alpha_H = C d\theta$ and $\lambda_\epsilon = 1$, hence
$\alpha_T = 0$ by \eqref{eq:21}, and this is the only case where this
might happen by \eqref{eq:22}.

\smallskip

Proposition \ref{prop:4:5} shows a large part of $\spec^*(d
*_\varepsilon)$ is symmetric with respect to $\frac{1}{2}$ and is
parametrised trough \eqref{eq:20} by the spectrum $\{ k + \epsilon
n^2\}$ of the non-negative elliptic Laplacian $L_{\varepsilon,H} =
\Delta_H - \varepsilon T^2 $ acting on $\Omega^2 H$, or equivalently
by the spectrum of
\begin{equation}
  \label{eq:23bis}
  \Delta_\varepsilon = \Delta_H - \varepsilon T^2
\end{equation}
acting on functions.  However there are ``holes'' in this symmetry
corresponding to the eigenvalues $\lambda_\varepsilon^- = -\varepsilon
k$ when $k = |n|$, for in this case $\alpha_T$ is not uniquely
determined by $\alpha_H$ so that we will have to treat these on a
separate footing.  This means that in the case $\lambda_\epsilon =
\lambda_\epsilon^-$, we have to remove from the parameter space the
horizontal forms $\alpha_H$ in
\begin{equation}
  \label{eq:24}
  \mathcal{H}^0 = \ker (\Delta_H^2 + T^2).
\end{equation}
This space has a simple description using the complex structure $J$
and the associated splitting $\Omega^1 H \otimes \C = \Omega^{1,0} H
\oplus \Omega^{0,1} H$. We recall that the component $d_H^{0,1}$ of
$d_H$ from functions to $ \Omega^{0,1} H$ is called the $\overline
\partial_b$ operator, and its kernel is the space of CR functions.

\begin{prop}
  \label{prop:4:6}
  The space $*_H \mathcal{H}^0$ is the space of pluri-CR functions,
  \emph{i.e.}  real parts of CR functions.
\end{prop}

\begin{proof}
  Consider the Kohn Laplacians $\overline{\square}_b = \overline
  \partial_b^* \overline \partial_b $ and $\Box_b = \partial_b^*
  \partial_b$ acting on functions.  Following, say, \cite[Theorem
  2.3]{Lee86}, one has in dimension $3$
  \begin{equation}
    \label{eq:25}
    \Delta_H = \overline{\square}_b + \square_b \ \quad \mathrm{and}
    \quad iT = \overline{\square}_b - \square_b.
  \end{equation}
  Since $T$ commutes with everything here one gets
  $$
  \Delta_H^2 + T^2 = 4 \overline{\square}_b \square_b = 4 \square_b
  \overline{\square}_b.
  $$
  If $f$ is a real function in $\mathcal{H}^0$ then $g=\square_b f$
  is CR since its image by $\overline{\square}_b$ is zero, and is in
  the image of $\Delta_H$ since its integral vanishes. Hence
  $$
  \Delta_H f = \overline{\square}_b f + \square_b f = \bar g + g =
  2 \Re g,
  $$
  and $f=2\Re h$ with $h = \Delta_H^{-1} g$ is a CR function as
  needed.
\end{proof}

\smallskip

We now study the missing case $\lambda_\varepsilon^2 = \varepsilon^2
n^2$.  We first recall that complex vertical forms $\Omega^* V \otimes
\C \simeq \theta \land \Omega^* H \otimes \C$ also have a natural
bigrading inherited from $J$ on $H$, independently from $\theta$.
Of particular interest here is the

\begin{defi}
  \label{defn:4:7}
  The bundle $K_M \simeq \theta \land \Omega^{1,0} H $ of $2$-forms
  vanishing on $H^{0,1}$ is called the canonical CR bundle.  We denote
  by $\mathcal{H}^{2,0}$ its subspace of closed sections, also called
  holomorphic $(2,0)$-forms, and $\mathcal{H}^2_+$ the real part of
  $\mathcal{H}^{2,0}$.
\end{defi}

When the CR manifold $M$ can be locally embedded in a $4$-dimensional
complex manifold $N$, $K_M$ is the restriction to $M$ of the canonical
bundle $K_N= \Omega^{2,0} N$ of $N$, and holomorphic forms are local
restrictions of holomorphic $(2,0)$-forms in $N$, see \cite{Lee86} for
instance.  This explains the notation in the previous definition, as
$\mathcal{H}^{2,0}$ (resp. $\mathcal{H}^2_+$) is related to the space
of holomorphic $(2,0)$-forms in the usual sense on $N$ (resp. to the
space of self-dual $2$-forms, orthogonal to the Kähler form). Note
that this is indeed the case for our CR-Seifert manifolds for one can
take $N = M × \R$ with the extension of $J$ considered above.

\smallskip

We now show that the remaining spectrum of $d *_\varepsilon$ is
entirely given by holomorphic forms.

\begin{prop}
  \label{prop:4:8}
  A $2$-form $\alpha \in \Omega^2_d M$ satisfies
  \begin{equation}
  \label{eq:17ter}
  (d *_\varepsilon )\alpha = \lambda_\varepsilon^- \alpha \ , \quad
  K\alpha = k\alpha \quad \text{and} \quad T^2\alpha  = -n^2\alpha 
  \end{equation}
  with $(\lambda_\epsilon^-)^2 = \epsilon^2n^2$ \upn{(}\emph{i.e.}
  $k=|n|$\upn{)} if and only if $\alpha_H=0$ and $\alpha = \theta
  \land \alpha_T $ belongs to $\mathcal{H}^2_+$.
\end{prop}

\begin{proof}
  Let $\alpha = \theta \land \alpha_T + \alpha_H$ in $\Omega^2_d M$
  be an eigenfunction of $d *_\varepsilon$ satisfying \eqref{eq:17ter}
  and $\lambda_\varepsilon^2 = \varepsilon^2 n^2$. By \eqref{eq:19}
  one has also $ \lambda_\varepsilon =-\varepsilon k$.  Since $(T*_H)^
  2= -T^ 2 = n^ 2 = k^2$ on $\Omega^1 H$, one can split
  $$
  \alpha_T = \alpha_T^+ + \alpha_T^- \quad \mathrm{with} \quad
  (T*_H) \alpha_T^ ± = ± k \alpha_T^ ±.
  $$
  Then \eqref{eq:21} is equivalent to
  \begin{equation}
    \label{eq:26}
    2 \varepsilon\, k \, \alpha_T^+ =   d_H *_H \alpha_H. 
  \end{equation}
  Moreover $K \alpha= k \alpha$ gives $(\Delta_H + T*_H ) \alpha_T = k
  \alpha_T$, which implies $\Delta_H \alpha_T^+ =0$ since $[\Delta_H,
  T*_H]=0$ on $\Omega^ 1H$. Therefore $\alpha_T^+$ lives in $\ker
  \delta_H$ leading by \eqref{eq:26} to $\Delta_H*_H \alpha_H = 0$,
  hence to $\alpha_H = C d\theta$ and $k=n=0$. If $C\not=0$, this
  implies by the eigenvalue identity \eqref{eq:19} that either
  $\lambda_\epsilon = \lambda_\epsilon^- = 0$, which is impossible
  since we consider the non-zero spectrum, or to $\lambda_\epsilon =
  \lambda_\epsilon^+ = 1$ which is impossible, too, because one would
  have $(\lambda_\epsilon)^2\not= \epsilon^2n^2$. We get then that in
  any case considered in the present proof, $\alpha_H=0$, so that
  $\alpha$ is a vertical form.
  
  Now \eqref{eq:22} reads $\delta_H \alpha_T=0$, or equivalently
  $$
  d_H (\theta \land J\alpha_T) = 0.
  $$
  Recall now that $\alpha$ belongs to $\Omega_d^2M$, hence is
  closed.  The $(1,0)$-part of $\alpha_T$ is then closed and $\theta
  \land \alpha_T$ lives in $\mathcal{H}^2_+$ as needed.
  
  Conversely, $\mathcal{H}^2_+$ is preserved by $J$ and $T$. Thus it
  can be split in eigenspaces of $T *_H = -J T = k$, on which $d
  *_\varepsilon = k$ by definition, see \eqref{eq:12}.
\end{proof}

We now summarize our spectral study of $d *_\varepsilon$ in relation
to the diabatic limit $\varepsilon \to 0$.

\begin{cor}
  \label{cor:spectre}
  The spectrum of $d *_\varepsilon$ splits into the following
  families:
  
  \begin{enumerate}
   
  \item A converging part $\Lambda_\varepsilon^+$, converging to $1$
    and parametrised by the whole spectrum of $\Delta_\varepsilon =
    \Delta_H - \varepsilon T^2$ \upn{(}acting on functions\upn{)} by
    the formula
    $$
    \Lambda_\varepsilon^+ = \spec \left( \frac{1 + \sqrt{1+ 4
          \varepsilon \Delta_\varepsilon}}{2} \right).
    $$
    
  \item A collapsing part, converging to $0$, itself divided into two
    families:

     \begin{enumerate}
    
     \item the first one $\Lambda_\varepsilon^-$, nearly symmetric to
       $\Lambda_\varepsilon^+$:
       $$
       \Lambda_\varepsilon^- = \spec \left(\frac{1 - \sqrt{1+ 4
             \varepsilon \Delta_\varepsilon}}{2} \right),
       $$
       but $\Delta_\varepsilon$ has here to be restricted to the
       orthogonal of the space of pluri-CR functions $\mathcal{H}^0$.
  
     \item the spectrum $\Lambda_\varepsilon^0$ of $\varepsilon T *_H
       = -\varepsilon JT$ acting on $\mathcal{H}^2_+$, the real
       parts of holomorphic forms in the canonical CR bundle.
      
      \end{enumerate}
  \end{enumerate}
\end{cor}

The signs of the eigenvalues in the first two families are clear.
About the third one, we can notice:

\begin{prop}\label{prop:H2+}
  Up to some finite dimensional space, $d*_\varepsilon$ is
  \emph{positive} on $\mathcal{H}^2_+$.
\end{prop}

\begin{proof}
  Recall that $d*_\varepsilon = -JT$ on $\mathcal{H}^2_+$. Consider
  then the splitting of the Tanaka-Webster connection $\nabla_H =
  \nabla_{1,0} + \nabla_{0,1}$ on $H \otimes \C$. Then on $K_M =
  \theta \land \Omega^{1,0} H$ one has in dimension $3$,
\begin{displaymath}
  R = \nabla_{0,1}^* \nabla_{0,1} - \nabla_{1,0}^* \nabla_{1,0} - i\nabla_T. 
\end{displaymath}
On holomorphic forms $\mathcal{H}^{2,0}$ in $K_M$, the Lie derivative
in $T$ equals $\nabla_T$ and the previous equation reduces to
$$
-iT = R + \nabla_{1,0}^* \nabla_{1,0},
$$
which implies that $-(iT + R)$ is a non-negative operator. As the
spectrum of $d*_\varepsilon$ (on closed forms) is discrete and without
accumulation points, there is only a finite dimensional space of
eigenvectors with nonpositive eigenvalues.
\end{proof}

In order to get more symmetry in the spectral decomposition of $d
*_\varepsilon$, one can fill in the holes in $\Lambda_\varepsilon^-$
by adding $\Delta_\varepsilon$ on $\mathcal{H}_0$. As already
discussed, this corresponds to adding the cases $k= |n|$ and
$\lambda_\varepsilon = -\varepsilon k \not= 0$. Given $k$, the
multiplicity of each added \emph{virtual} eigenvalue $-\varepsilon k$
is equal to $2 h_0(k)$ by Proposition \ref{prop:4:6}, where we have
denoted
$$
h_0(k) = \mathrm{dim}_\C\,\bigl\{ \text{CR functions } f \text{
  such that } iT f= - kf \bigr\}.
$$
Observe that by \eqref{eq:25}, $h_0(k)=0$ if $k<0$.  In the same
spirit, the holomorphic part $\Lambda_\varepsilon^0$ above consists in
$\{ \varepsilon k \mid k \in \Z^*\}$, with multiplicity $2 h_2(k)$
given by
$$
h_2(k) = \mathrm{dim}_\C\,\bigl\{ \text{holomorphic }
(2,0)\text{-forms } \alpha \in \mathcal{H}^{2,0} \text{ such that } iT
\alpha= - k \alpha \bigr\}.
$$
Considering the positive operators
$$
Q_\varepsilon^± =  \frac{± 1 + \sqrt{1+ 4 \varepsilon
      \Delta_\varepsilon}}{2 \varepsilon}  \ ,
$$
leads to the more suggestive decomposition:
\begin{equation} \label{eq:27}
  \spec^* \left( \frac{d *_{\varepsilon}}{\varepsilon} \right) = 
  ± \spec^* \left(Q_\varepsilon^±\right) \ \cup \ 2×
  \spec^* \left( - iT_{\mid \mathcal{H}^{2,0}} \right)\ 
  \setminus \ 2 × \spec^* \bigl( iT_{\mid \ker \overline\partial_b}
  \bigr). 
\end{equation}

This formula shows that the virtual spectrum of $d*_\epsilon$ consists
in a two completely different parts: a (nearly) symmetric part to
$1/2$, that varies with $\epsilon$, and a \emph{constant} holomorphic
part. We will see in Lemma \ref{lem:zetaQ_pm} that the symmetric part
always contributes to $1$ in the renormalized $\eta$-invariant
$\eta_0$ when torsion vanishes. Hence the computation of $\eta_0$
finally reduces to counting holomorphic objects, as will be done in
section \ref{sec:decompostion-sp-d}. This phenomenon has already been
observed on a smooth base in \cite{Zhang} and over orbifolds, in the
adiabatic context and constant curvature, in \cite{Nicolaescu2}.

\bigskip

\section{The spectrum of $D*$ and comparison of the 
  $\eta$-invariants} \label{sec:decompostion-sp-d}

Our goal is now to relate our description of the spectrum of
$P_\varepsilon = \varepsilon^{-1}d *_\varepsilon$ to the spectrum of
the middle operator of the contact complex $D*$. We already know (see
the discussion at the end of section \ref{sec:4:1}) that the bounded
spectrum of $P_\varepsilon$ converges towards that of $D*$ in the
diabatic limit \cite{Rumin00}.
Therefore from Corollary \ref{cor:spectre} the non-zero spectrum of
$D*$ has to split as follows
\begin{equation}
  \label{eq:28}
  \spec^*(D*) = \spec^*(-\Delta_H \ \text{on} \ (\mathcal{H}^{0})^\bot)
  \ \cup \ \spec^*(-J T \ \text{on} \ \mathcal{H}^2_+) 
\end{equation}
(note the lack of uniformity already noted in the introduction in the
convergence of $\Lambda_\epsilon^-$ when $\epsilon\to 0$, as each
eigenvalue $\mu$ in the spectrum of $\Delta_H$ is approached at a
speed approximately $\epsilon\mu$).  This is enough to compare the
needed $\eta$-invariant to $\eta_0$ and conclude (see \eqref{eq:31}
below and the discussion following it), but we would like first to
spend a few lines to reinterpret this more precisely in the CR-Seifert
context.

\smallskip
\subsection*{The spectrum of $D*$ from the CR viewpoint}

First of all, the second spectral family of eigenvalues in
\eqref{eq:28} is clearly embedded in $\spec^*(D*)$, as \eqref{eq:13}
shows that $D* = -TJ$ on $\mathcal{H}^2_+$. To understand where the
first one comes from, we consider the following operator
$$
Q = d_H J : \ker d_H \subset \Omega^2 V \longrightarrow \Omega^3 M
.
$$
By definition $\mathcal{H}^2_+ = \ker Q$. We also remark that
$$
(Q^*) *_M = (\Pi_{\ker d_H}J \delta_H) *_M = -*_M (\Pi_{\ker
  d_H}J d_H)
$$
so that $\ker Q^* = *_M \mathcal{H}^0$ and $\overline{\im Q} = *_M
\left( \mathcal{H}^0\right)^\bot $.  To complete the landscape, we of
course define $\mathcal{H}^2_- = \overline{\im Q^*}$, so that
\begin{equation}
  \label{eq:29}
  \ker d_H \cap \Omega^2 V= \ker Q \oplus \overline{\im Q^*} =
  \mathcal{H}^2_+ \oplus \mathcal{H}^2_-.
\end{equation}
Then in vanishing Webster torsion, one has by \eqref{eq:13} that
\begin{equation}
  \label{eq:30}
  \begin{split}
    Q (D*) & = d_H J (-TJ - (d_H *_H)^2)  = T d_H  + (d_H *_H)^3 \\
    & = - \Delta_H Q ,
  \end{split}
\end{equation} 
on $\ker d_H \subset\Omega^2 V$, where $\Delta_H = d_H \delta_H$ is
the contact Laplacian on $\Omega^3 M$, conjugate to $\Delta_H$ on
functions through $*_M$. This shows that $D*$ is conjugate to
$-\Delta_H$ on $*_M \left( \mathcal{H}^0\right)^\bot$ by $Q$, and that
$D*$ preserves the splitting \eqref{eq:29}. We therefore recover the
decomposition of $\spec(D*)$ in two families \eqref{eq:28}, but now
entirely seen within $\Omega^2 V$ :
\begin{equation}
  \label{eq:29bis}
  \spec^*(D*)  = \spec^*(D*_{\mid \mathcal{H}^2_- = \overline{\im Q^*}})
  \ \cup \ \spec^*(D*_{\mid  \mathcal{H}^2_+ = \ker Q}) .
\end{equation}
where by \eqref{eq:30}, $D*$ is conjugate to $-\Delta_H$ on $*_H(
\mathcal{H}^0)^\bot$ by $Q$.

\smallskip

The space $\mathcal{H}^2_-$ is actually a CR invariant, as is
$\mathcal{H}^2_+$. Indeed $\Delta_H$ is surjective on $\Omega^3 M$ up
to ``constant'' $3$-forms $C \theta\land d\theta$; as $Q^*$ is zero
on these,
\begin{align*}
  \mathcal{H}^2_- & = \overline{\im Q^*} = \overline{\im Q^*\Delta_H} \\
  &= \overline{\im D* J \delta_H} \ ,\quad \text{by} \quad
  \eqref{eq:30}, \\
  & = \overline{\im D J d_H}.
\end{align*}

\smallskip 

We now have two splittings of $\Omega^2 V\cap \im D$ : the spectral
one
$$
\im D = E^+ \oplus E^-
$$
in the positive and negative eigenspaces of $D*$, and the CR
invariant one given by
$$
\im D = (\mathcal{H}^2_+ \cap \im D) \oplus \mathcal{H}^2_-\ .
$$
It follows from Prop. \ref{prop:H2+}, \eqref{eq:28} and
\eqref{eq:29} that, on a CR-Seifert manifold, the pair $(E^+, E^-)$ is
in Fredholm position with respect to $(\mathcal{H}^2_+,
\mathcal{H}^2_-)$. More precisely,
$$
\mathcal{H}^2_+ = E^+ \oplus V \oplus H^2(M, \R) \quad \text{and} \quad E^- =
\mathcal{H}^2_- \oplus V
$$
with the finite dimensional space $V = \mathcal{H}^2_+ \cap E^-$.
This enlightens the CR meaning of the spectral asymmetry of $D*$ we
are studying here.

Observe however that if the formal definitions of $\mathcal{H}^2_±$
make sense on any $3$-dim CR manifold, their use is highly problematic
in general. For instance $\mathcal{H}^2_+$ may be empty if $M$ does
not bound a Stein manifold, while $E^+$ and $E^-$ still exist and keep
their nice analytic features by hypoellipticity of $D*$ on $\im D$.
The previous Fredholm picture then definitely breaks down. Anyway, from
the pseudodifferential viewpoint, the projection on $E^+$ is a natural
quantization of the real part of the Szegö projector on holomorphic
$(2,0)$-forms, as seen at the Heisenberg symbolic level, see e.g
\cite[Chap 4]{Be-Gr} for more details on this notion.

\smallskip

We now come back to the comparison between the Riemannian and contact
spectra. In \eqref{eq:28}, we can proceed as in \eqref{eq:27} by
``filling the holes'' in the spectrum of $-\Delta_H$ on
$\mathcal{H}^0$. From \eqref{eq:25} we still have $\Delta_H= -iT$ on
CR functions, and this leads to the following decomposition:
\begin{equation}
  \label{eq:31}
  \spec^*(D*) \ = \ \spec^*(-\Delta_H) 
  \ \cup \ 2× \spec^*(- iT_{\mid \mathcal{H}^{2,0}} ) 
  \ \setminus \ 2× \spec^* (iT_{\mid \ker \overline\partial_b}).  
\end{equation}

\begin{rem} \label{rem:symetrie_D}
  In a slightly more tricky way, one can add $\spec^*(\Delta_H)$ to
  both sides of \eqref{eq:31}: the operator $\Delta_H$ on functions is
  conjugate to $\Delta_H = d_H \delta_H$ on $\Omega^3 M$ and, wedging
  by $\theta$, to $\delta_Hd_H$ on $\Omega^2 V$. The spectrum of the
  contact Laplacian
  $$
  \Delta_2 = D* + \delta_H d_H \quad \textrm{on } \Omega^2 V
  $$
  (see section \ref{sec:eta+} for more on this one) appears then in
  a very symmetric manner, namely
  \begin{equation}
  \label{eq:32}
  \begin{split}
    \spec^* (\Delta_2) & = \spec^* (D*) \ \cup \ \spec^*(\Delta_H) \\
    &  = \spec^*(\Delta_H) \ \cup \ \spec^*(- \Delta_H) \\
    & \quad \quad \bigcup 2× \spec^*(- iT \mid \mathcal{H}^{2,0}
    ) \setminus 2× \spec^* (iT \mid \ker \overline\partial_b) ,.
  \end{split}
  \end{equation}
  This spectral symmetry can also be seen directly. Equation
  \eqref{eq:13} yields
  $$
  \Delta_2 = T*_H - d_H \delta_H + \delta_H d_H = T*_H + \, P
  $$
  on $\Omega^2 V = \theta \land \Omega^1 H$. As $[*_H, T*_H] = 0$
  while $*_H P = -P *_H$, $\Delta_2 (*_H P) = - (*_H P) \Delta_2$ and
  $\spec(\Delta_2)$ is symmetric except maybe on $\ker P$, where
  $\Delta_2= T*_H = -TJ$. It is then easily seen that the kernel
  splits into
  $$
  (\ker P)^{2,0} = \mathcal{H}^{2,0} \oplus
  {\overline\partial_b}^{-1}(*_M \ker \overline\partial_b) ,
  $$
  yielding \eqref{eq:32}.
\end{rem}

\begin{rem}
  Let us mention that this decomposition and the spectral symmetry of
  $\Delta_2$ also hold on contact manifolds of any dimension, in
  vanishing Tanaka-Webster torsion, see \cite[Prop.  8]{Rumin94}.
  This leads to the same kind of formulae as \eqref{eq:32}, with a
  ``residual spectrum'' given by sum of $\eta$-functions counting
  holomorphic objects.
\end{rem}

\subsection*{Comparison of contact and Riemannian eta invariants}

Comparing the spectrum of $P_\epsilon$ given by \eqref{eq:27} with
that of $D*$ in \eqref{eq:31} yields an immediate relation between
their $\eta$-functions, up to combinations of $\zeta$-functions of
positive operators:

\begin{prop} \label{prop:etaD_P_epsilon}
  On a CR-Seifert manifold,
  \begin{equation}
  \label{eq:33}
  \eta(P_\varepsilon) - \eta (D*) = \zeta(\Delta_H)
   + \zeta(Q^+_\varepsilon) - \zeta(Q^-_\varepsilon),
\end{equation}
where $\dsp Q_\varepsilon^± = \frac{1}{2\varepsilon}( ± 1 +
\sqrt{1 + 4 \varepsilon \Delta_\varepsilon})$, and
$\Delta_\varepsilon= \Delta_H - \varepsilon T^2$ on functions.
\end{prop}

\smallskip

Following Definition \ref{defi:eta_0}, the renormalized $\eta$ invariant
$\eta_0(M, \theta)$ is the constant term in the development of
$\eta(P_\varepsilon) (0) = \eta(M, g_\varepsilon)$ in powers of
$\epsilon$. It is then immediately extracted from \eqref{eq:33} as
follows:
 \begin{equation}
  \label{eq:34bis}
  \eta_0(M,\theta) = \eta (D*)(0) + \zeta(\Delta_H)(0) + \zeta_0(Q) ,
\end{equation}
where $\zeta_0(Q)$ is the constant term in the development in powers
of $\epsilon$
\begin{equation}\label{eq:34ter}
 \zeta(Q^+_\varepsilon)(0) - \zeta(Q^-_\varepsilon)(0) = \sum_{i=-2}^{2}
\zeta_i(Q) \,\epsilon^i ,
\end{equation}
which we already know to exist by \eqref{eq:2:6} and \eqref{eq:33},
since it is the same as that of $\eta(P_\epsilon)$ except for the
constant term. Moreover, it turns out that $\zeta_0(Q)$ can be
evaluated without too much harm on \emph{arbitrary} CR manifolds of
dimension $3$.

\begin{lem}
  \label{lem:zetaQ_pm}
  On any $3$-dimensional CR manifold,
  $$
  \zeta(Q_\varepsilon^+)(0) = - \zeta(Q_\varepsilon^-)(0),
  $$
  and
  $$
  \zeta_0(Q) = \frac{1}{24\pi^2} \int_M |\tau|^2 \theta \land d
  \theta.
  $$
  where $\tau = -\frac{1}{2} J \mathcal{L}_T J$ is the
  Tanaka-Webster torsion.
\end{lem}

\begin{proof}
  In view of
  $$
  2 \varepsilon \, Q_\varepsilon^± = ± 1 + \sqrt{1 + 4 \varepsilon
    \Delta_\varepsilon},
  $$
  we consider for $\lambda \geq -1$ the family of operators
  $$
  Q (\lambda) = \lambda + \sqrt{1 + 4 \varepsilon
    \Delta_\varepsilon},
  $$
  where actually
  $$
  \varepsilon \Delta_\varepsilon = \varepsilon \Delta_H -
  \varepsilon^2 T^2 = \Delta_{g_\varepsilon }
  $$
  is the standard Laplacian on functions for the rescaled metric
  $g_\varepsilon= \varepsilon^{-2} \theta^2 + \varepsilon^{-1}
  \gamma_H$ we use here.
  
  Seeley's classical results \cite{Seeley} infer that $Q(\lambda)$ is
  a smooth family of positive elliptic pseudo-differential operators
  of order $1$, and that their $\zeta$-functions
  $$
  P(\lambda)(s) := \zeta (\lambda + \sqrt{1 + 4
    \Delta_{g_\varepsilon}})(s)
  $$
  are meromorphic with possibly simple poles at $s=1$, $2$ and $3$.
  According to \cite[Prop. 2.9]{aps3} or \cite[Lemma 1.10.2]{Gilkey}
  one can differentiate $P(\lambda)(s)$ with respect to $\lambda$ to
  get
  $$
  \frac{d }{d\lambda} P(\lambda)(s) = - s P(\lambda)(s+1).
  $$
  Therefore $\dsp \frac{d^4 }{d\lambda^4} P(\lambda)(0) = 0$ since
  $P(\lambda) $ is regular at $s=4$, and $P(\lambda)(0)$ is a
  polynomial of degree $3$ in $\lambda$:
  \begin{equation}
    \label{eq:34}
    P(\lambda) = \zeta((1+ 4 \Delta_{g_\varepsilon})^{1/2})(0) -
    \lambda R_1 + \lambda^2 
    \frac{R_2}{2} - \lambda^3 \frac{R_3}{3},
  \end{equation}
  where $R_0=\zeta(\sqrt{1+4\Delta_{g_\epsilon}})(0)$ and $R_n$ for
  $n>0$ stands for the residue at $s=n$ of
  $$
  \zeta(\sqrt{1+ 4 \Delta_{g_\varepsilon}})(s) = \zeta(1+ 4
  \Delta_{g_\varepsilon})(s/2).
  $$
  Actually these residues are related to the development of the
  heat kernel of $\Delta_{g_\epsilon}$ on functions in a simple way.
  Let
  \begin{displaymath}
    \tr(\exp^{-t \Delta_{g_{\varepsilon}}}) \stackrel{t
    \to 0^+}{\sim}  
    \frac{a_0(g_\varepsilon)}{t^{3/2}} + 
    \frac{a_2(g_\varepsilon)}{t^{1/2}} + \cdots .
  \end{displaymath}
  According to \cite[Theorem 4.8.18d]{Gilkey}, the constants are computed
  in terms of the volume and the Riemannian scalar curvature of
  $g_\varepsilon$ as:
  \begin{equation}
    \label{eq:35}
    a_0(g_\varepsilon) = \frac{\mathrm{Vol}(M,g_\varepsilon)}{(4\pi)^{3/2}} 
    \quad \text{and} \quad
    a_2(g_\varepsilon) = \frac{1}{6 (4\pi)^{3/2}}\,\int_M
    \Scal(g_\varepsilon) d \vol_{g_\varepsilon} \ .
  \end{equation}
  This yields
  $$
  \tr(\exp^{-t(1 + 4\Delta_{g_\varepsilon})}) = e^{-t} \tr(\exp^{-4
    t\Delta_{g_\varepsilon}}) \sim \frac{a_0(g_\varepsilon) }{8
    t^{3/2}} + \frac{4 a_2(g_\varepsilon) -a_0(g_\varepsilon)
  }{8t^{1/2}} + \cdots ,
  $$
  and by Mellin's transform \cite[Lemma 1.10.1]{Gilkey},
  \begin{equation*}
    \Gamma(s/2)\, \zeta(1+\Delta_{g_\varepsilon})(s/2) 
    \ = \ \frac{a_0(g_\varepsilon) }{4(s - 3)} \ + \ \frac{ 4
      a_2(g_\varepsilon) -a_0(g_\varepsilon) }{4(s - 1)} \ + \ h(s),
  \end{equation*}
  with $h$ holomorphic for $\Re s > -1$. Hence
  $$
  \zeta((1+ 4 \Delta_{g_\varepsilon})^{1/2})(0) = 0 \,
  $$
  as this is the only way to cancel the simple pole of the Gamma
  function at $s=0$, and
  $$
  R_2 = 0,
  $$
  (because the Gamma function does not vanish at $s=2$ and the
  r.h.s. has no pole at this point) so that $P(\lambda)$ is an odd
  polynomial.  This gives $P(1) = -P(-1)$ or, equivalently,
  $$
  \zeta(Q_\varepsilon^+)(0) = -\zeta(Q_\varepsilon^-)(0)
  $$
  as announced.  Moreover one has
  $$
  R_1 = \frac{4\,a_2(g_\varepsilon ) - a_0(g_\varepsilon)}{4\sqrt
    \pi} \quad \mathrm{and} \quad R_3 =
  \frac{a_0(g_\varepsilon)}{2\sqrt{\pi}},
  $$
  and thus by \eqref{eq:34} and \eqref{eq:35}
  \begin{equation}
    \label{eq:36} 
    \begin{split}
      \zeta(Q_\varepsilon^+)(0) & = -R_1 - R_3/3 \\
      & = \frac{1}{\sqrt \pi}
      (\frac{a_0(g_\varepsilon)}{12} - a_2(g_\varepsilon)) \\
      & = \frac{1}{48 \pi^2 \varepsilon^2} \,\bigl(
      \frac{1}{2}\,\int_M \theta \land d \theta \, - \int_M
      \Scal(g_\varepsilon) \,\theta \land d\theta \bigr).
    \end{split}
   \end{equation}
   The Riemannian curvature of $g_\varepsilon$ can be developed in
   powers of $\varepsilon$ using the links between Tanaka-Webster and
   Levi-Civita connections underlined in \eqref{2:6bis}. According to
   \emph{e.g.}  \cite[p 318]{Rumin94}, one finds in dimension $3$ that
   $$
   \Scal(g_\varepsilon) = -\frac{1}{2} + 2 \varepsilon R -
   \varepsilon^2 \, |\tau|^2,
   $$
   where $R$ and $\tau$ are Tanaka-Webster curvature and torsion.
   The constant term in the full development of
   $\zeta(Q_\varepsilon^+)$ is then necessarily equal to the integral
   of $\frac{1}{48\pi^2}\, |\tau|^2$ on $M$.
\end{proof}

\begin{rem}
  \label{rem:conjecture}
  According to \eqref{eq:27}, $Q_\varepsilon^+$ describes the non
  collapsing spectrum of $d *_\varepsilon$, on Seifert-CR manifolds.
  We have seen that this spectrum only contributes by a local
  expression $\zeta(Q_\varepsilon^+)(0) $ to $\eta(d *_\varepsilon)$.
  We expect this to hold in the general case. Indeed on any CR
  manifold, the non-collapsing spectrum is always strictly positive,
  since it converges to $1$ and $d*_\varepsilon$ has no spectral flow.
  It therefore always contributes through a \emph{zeta} function,
  whose value at $0$ is local for a wide class of operators.
\end{rem}

\smallskip

\subsection*{A computation of $\eta_0$.}\label{rem:calculNicolaescu2}

The previous Lemma \ref{lem:zetaQ_pm}, together with the spectral decomposition
\eqref{eq:27}, leads to a general computation of the renormalized
$\eta$-invariant on all CR-Seifert manifolds, including the still
missing exceptional cases of section \ref{sec:Sfrt}. Indeed, one has
$$
\zeta^*(Q_\varepsilon^+) - \zeta^*(Q_\varepsilon^-) =
\zeta(Q_\varepsilon^+) - \zeta(Q_\varepsilon^-) +1,
$$
since $0$ belongs to $\spec(Q_\varepsilon^-)$ with multiplicity $1$
(corresponding to the constant functions). It follows then from
\eqref{eq:27} that
\begin{equation}
  \label{eq:37}
  \eta_0(d *) = \eta_{\mathrm{ad}}(d*) = 1 + 2 \, \bigl( \eta(-iT_{\mid
  \mathcal{H}^{2,0}})(0)
  - \eta(iT_{\mid \ker \overline\partial_b})(0) \bigr).
\end{equation}

These holomorphic counting functions can be nicely expressed as
dimensions of spaces of sections on adequate orbifold line bundles
over the basis orbifold Riemann surface, which in turn are easily
computed with the help of Riemann-Roch-Kawasaki's theorem
\cite{kawasaki}.  Note that this has already been observed in the
adiabatic setting and constant curvature by L. Nicolaescu in
\cite[Sec. 1]{Nicolaescu2}. We give below only a short description of
the computation, and refer to \cite{Nicolaescu2} for more details.

\smallskip

Following section \ref{sec:Sfrt}, the CR-Seifert manifold $M$ may be
seen as the unit circle bundle of some orbifold line bundle $L$ over
$\Sigma$, with singular data $(\alpha_i, \rho_i, \beta_i)$ at points
$m_i \in \Sigma$. 
Let $K_\Sigma = \Lambda^{1,0} T^* \Sigma$ denotes the orbifold
canonical bundle of $\Sigma$. Now, given a Fourier component $i T = n
\in \Z$, the space of CR functions $f$ such that $f \circ\varphi_t =
e^{-int} f$ are interpreted as the space of holomorphic sections of
$L^n$, and we denote by $h_0(L^n)$ its dimension. Moreover the space
of holomorphic forms $\sigma$ in the canonical CR bundle $K_M \simeq
\theta \land K_\Sigma \otimes L$ such that $-iT \sigma = n \sigma$ may
be seen as the space of holomorphic sections of $K_\Sigma \otimes
L^n$, \emph{i.e.}  $(1,0)$-holomorphic forms in $L^n$. Let $h_1(L^n)$
denotes its dimension. Hence we get
\begin{equation}
    \label{eq:38}
    \begin{aligned}
      \eta(-iT_{\mid \mathcal{H}^{2,0}})(s) - \eta(iT_{\mid \ker
        \overline\partial_b})(s) & = - \sum_{n \in \Z^*}
      \mathrm{sgn}(n)
      \frac{h_0(L^n) - h_1(L^n)}{|n|^s} \\
      & = \sum_{n \in \Z^*} \mathrm{sgn}(n)
      \frac{\chi_{\overline\partial} (L^{-n})}{|n|^s}.
    \end{aligned}
\end{equation}
Following the method in \cite[Sec. 1]{Nicolaescu2}, this sum can be
computed explicitly using Riemann-Roch-Kawasaki theorem (extension of
the classical Riemann-Roch to the orbifold case) \cite{kawasaki}.
Using the (rational) orbifold Euler characteristic $\chi$ of the base
$\Sigma$ and the (rational) degree $d$ of $L$, it reads
\begin{equation}
  \label{eq:39}
  \chi_{\overline\partial}(L^{-n}) = \frac{\chi}{2} - n d + \sum_i
  \frac{1}{2}\left( 1 - \frac{1}{\alpha_i} \right) - \Bigl\{ \frac{-n
  \beta_i \rho'_i}{\alpha_i} \Bigr\},
\end{equation}
where $\{x\} = x - [x]$ denotes the fractional part of $x$, and
$\rho'_i$ is the inverse of $\rho_i$ mod.~$\alpha_i$. This purely
topological formula holds true, irrespective of the curvature value.
The result should then be the same in the constant and non-constant
curvature cases, so that Ouyang's formula \eqref{eq:2} for $\eta_0$
holds true on any CR-Seifert manifold.

To get explicitly the formula, one can argue as follows: the constant
terms in \eqref{eq:39} do not contribute to the sum \eqref{eq:38},
whereas
$$
\sum_{n \in \Z^*} -d  |n|^{-s+1}=  -2 d \,\zeta(s-1)
$$
has value $\frac{d}{6}$ at $s=0$. The Dedekind-Rademacher sums
$s(\alpha_i, 1, \beta_i \rho'_i) = s(\alpha_i, \rho_i, \beta_i)$
appear from the periodic orbifold contribution in \eqref{eq:39}, as in
Nicolaescu's work using \cite[Proposition 1.4]{Nicolaescu2}.
Inserting in \eqref{eq:37} leads to the desired expression.

\begin{rem}  
  This last computation shows that Theorem \ref{th:D*eta0} could have
  been proved in a quicker way on constant curvature CR-Seifert
  manifolds: applying the previous formulae and using the computation
  of $\zeta(\Delta_H)(0)$ given below leads to an expression for
  $\eta(D*)$ that can be compared directly to Ouyang's formula for
  $\eta_0$.  We have however omitted this proof since the links
  between $\eta(D*)$ and $\eta_0$ proved in this way would have
  appeared as the result of a possibly completely fortuitous or
  miraculous equality between explicitly known numerical expressions.
  On the contrary, our proof stresses the fact that the relation
  between $D*$ and $d*$ is deeply rooted in the nature of CR geometry
  and the diabatic limit. Moreover, it applies to \emph{the whole
    family} of CR-Seifert manifolds, irrespective of their curvature,
  and especially the exceptional cases that do not admit constant
  curvature contact forms.
\end{rem}

We now complete the comparison between $\eta_0$ and the contact
$\eta$-invariant $\eta(D*)$.
\begin{theo}
  \label{thm:etaD_nu_0}
  Let $M$ be a CR-Seifert manifold. Then,
  \begin{equation}
    \label{eq:40}
    \eta_0(M, \theta) = \eta (D*)(0) + \zeta(\Delta_H)(0) 
  \end{equation} 
  with
  \begin{equation}
    \label{eq:41}
    \zeta(\Delta_H)(0) = \frac{1}{512}\,\int_M R^2 \,\theta \land d
    \theta \ . 
  \end{equation}
\end{theo}
  
\begin{proof}
  From Proposition \ref{prop:etaD_P_epsilon} and Lemma
  \ref{lem:zetaQ_pm} it remains to compute $\zeta(\Delta_H)(0)$.  The
  development of the heat kernel $\exp^{- t\Delta_H}$ of the Kohn
  Laplacian $\Delta_H$ has been studied by Beals, Greiner and Stanton
  in \cite[Theorem 7.30]{Be-Gr-St}. On any CR manifold of dimension $3$,
  $$
  \tr(\exp^{-t \Delta_H}) \ \sim\ \sum_{n=0}^\infty t^{n-2} b_n(M,
  \theta) \quad \text{as} \quad t \to 0^+,
  $$
  where $b_n(M, \theta)$ are integrals over $M$ of polynomials of
  covariant derivatives of Tanaka-Webster curvature and torsion.
  Mellin's transform yields again
  $$
  \Gamma(s)\, \zeta(\Delta_H)(s) \ = \ \sum_{n \leq N} \frac{b_n(M,
    \theta)}{s - 2 +n} \ + \ h_N(s)
  $$
  with $h_N$ holomorphic for $\Re s > N -2$, and hence
  $$
  \zeta(\Delta_H)(0) = b_2(M, \theta).
  $$
  As $\zeta(\Delta_H)(0)$ stays unchanged when $\theta$ becomes $k
  \theta$, one must have $b_2(M, k\theta) = b_2(M, \theta)$, and the
  same argument as in Lemma \ref{lem:2:3} gives that
  $$
  b_2(M, \theta)\ = \ C_1\,\int_M R^2 \,\theta \land d \theta \ +
  \ C_2\,\int_M |\tau|^2\, \theta \land d \theta,
  $$
  for some constants $C_1$, $C_2$.
  
  Thanks to N. Stanton's work \cite{Stanton} it is possible to
  determine $C_1$ on the sphere $\cerc^3$. Indeed, let $L= 4 \Delta_H
  + R$ be the CR-conformal Laplacian on $\cerc^3$. Stanton states in
  \cite[Theorem 4.34]{Stanton} that for the contact form $\theta= i
  {\overline \partial} r = \frac{i}{2}(z^1 d\bar z^1 + z^2 d\bar z^2)$
  $$
  \tr(\exp^{-t L}) = \frac{\pi^2}{256t^2} \, + \, O
  \bigl(\frac{1}{t^2}\,\exp^{-\pi^2/4t} \bigr) \quad \text{as} \quad t
  \to 0^+.
  $$
  Now Tanaka-Webster curvature $R = 4$ here, so that the heat
  development of $\Delta_H$ is
\begin{align*}
  \tr(\exp^{-t \Delta_H}) \ = \ \exp^t \tr(\exp^{-t L/4}) \ =\ 
  \exp^t \frac{\pi^2}{16 t^2} \, + \,O
  \bigl(\frac{1}{t^2}\,\exp^{-\pi^2/4t} \bigr),
\end{align*}
and the constant term is $b_2(M, \theta) = \frac{\pi^2}{32}$. Hence
$$
\zeta(\Delta_H)(0) \ = \ \frac{\pi^2 }{32} \ = \ C_1
\,\int_{\cerc^3} R^2 \,\theta \land d \theta \ = \ 16 \pi^2 \, C_1
$$
yields $C_1 = \frac{1}{32× 16}$ on the sphere, hence on any
CR-Seifert manifold.
\end{proof}

Putting together this last result and Theorem \ref{eta0} leads to
Corollary \ref{maintheo2}.

\bigskip

\section{The contact and the modified contact $\eta$-invariants}
\label{sec:eta+}

We first begin by showing existence of the contact $\eta$-invariant in
dimension $3$.  It follows mostly the classical method of Chapter 1 of
\cite{Gilkey}, using pseudo-differential calculi developed on contact
manifolds.  As a consequence, we shall put below the emphasis mainly
on the steps where the hypoelliptic context introduces differences
with the well-known elliptic theory.

\begin{theo} \label{thm:eta(D*)}
  Let $(M,H,J)$ be a compact $3$-dimensional strictly pseudoconvex
  CR manifold endowed with a compatible contact form $\theta$ and the
  associated metric $g_1 = \theta^2 + d\theta( \cdot, J \cdot)$.
  Then the series
  $$
  \eta(D*)(s) = \tr^* (D* |D*|^{-(s+1)} ) = \sum_{\lambda_i \in
    \spec(D*)\setminus \{0\} } \frac{\lambda_i}{|\lambda_i |^{s+1}}
  $$
  converges absolutely for $\Re s > 2$, and has an meromorphic
  extension with possible simple poles at $s = 2 - n/2$ for $n \in
  \N$. Moreover $\eta(D*)(s)$ is regular at $s=0$; its value
  $\eta(D*)(0)$ is the \emph{contact $\eta$-invariant}.
\end{theo}

\begin{proof}
  From \cite{Rumin94} the two Laplacians
  $$
  \Delta_2\ = \ D* \, + \,\delta_H d_H \ \text{on}\ \Omega^2 V
  \quad \text{and}\quad \Delta_3= d_H \delta_H \ \text{on}\ \Omega^3
  M,
  $$
  are maximally hypoelliptic (be careful: $\Delta_3$ is
  nonnegative, but $\Delta_2$ is not, despite the notation). This
  means that they control two horizontal derivatives in $L^2$ norms
  (and one vertical derivative).  By the associated Sobolev
  embeddings, their resolvents are compact and their spectra are
  discrete. By orthogonality and conjugation, the non-zero spectrum of
  $\Delta_2$ splits into
  \begin{equation}
    \label{eq:42}
    \spec^*(\Delta_2) = \spec^* (D*) \cup \spec^* (\Delta_3)\  ,
  \end{equation}
  and $D*$ has discrete pure point spectrum with finite multiplicities
  on $\im D$. Sobolev embeddings also yields that  $(i+
  \Delta_2)^{-n}$, $(i +\Delta_3)^{-n}$ are trace
  class for $n$ large enough, hence the same for $(D*)^{-n}$. The
  series $\eta(D*)(s)$ is then well defined and holomorphic for $\Re
  s$ large.
  
  Getting more information on $\eta$ relies in the Riemannian
  (elliptic) case on the use of the classical pseudo-differential
  calculus for elliptic operators. Such a symbolic calculus has also
  been developed on contact manifold by Beals, Greiner and Stanton in
  \cite{Be-Gr, Be-Gr-St} or Taylor in \cite{Taylor}, a concise account
  may also be found in \cite{Getzler}.  The symbols of the
  hypoelliptic operators $\Delta_2$ and $\Delta_3$ are invertible in
  this calculus: this follows from \cite[Lemmas 5.18, 5.19]{Ju-Ka}, or
  else by observing that in dimension $3$ their principal symbols are
  sums of invertible Folland-Stein ones.
  
  The parameter calculus adapted to the Heisenberg setting developed
  in propositions 5.20 to 5.26 of \cite{Ju-Ka} yields
  pseudo-differential approximations $R(\lambda)$ of the resolvents $(
  (\Delta_2)^2 - \lambda)^{-1}$, when $\lambda \notin \R^+$. This uses
  the classical iteration process described in \cite[p.  51]{Gilkey}
  or \cite[Sec. 9.1]{Shubin} for instance, where the standard
  pseudo-differential symbolic product has to be replaced by the
  Heisenberg one, see \cite{Be-Gr-St, Getzler}. 
  The symbol of these $R(\lambda)$ are universal expressions involving
  the symbol of $(\Delta_2)^2 - \lambda$, its inverse, and tensorial
  expressions of the Webster-Tanaka curvature and its derivatives.
  
  Then, as explained in \cite[Sec. 1.7]{Gilkey}, $R(\lambda)$ can be
  used in place of $( (\Delta_2)^2 - \lambda)^{-1}$ in the contour
  integral
  $$
  \Delta_2 \exp^{-t(\Delta_2)^2} \ = \ \frac{1}{2i\pi}\,
  \int_\gamma \exp^{-t\lambda} \Delta_2 ( \Delta_2^2 - \lambda)^{-1}\,
  d\lambda,
  $$
  with $\gamma \subset \C \setminus \R^+$ the correctly oriented
  boundary of the cone $\{\Im \lambda \leq \Re\lambda + 1\}$, in order
  to get good approximations of $\Delta_2 \exp^{-t(\Delta_2)^2}$ when
  $t$ goes to $0$. Following Lemma 1.7.7 of \cite{Gilkey}, homogeneity
  arguments then easily lead to the asymptotic development of
  $\Tr(\Delta_2 \exp^{-t\Delta_2^2})$ when $t \to 0^+$. Namely,
  \begin{equation}
    \label{eq:43}
    \Tr(\Delta_2 \exp^{-t\Delta_2^2}) \ \sim \ \sum_{n=0}^\infty t^{(n-6)/4} 
  R_n(M, \theta),
  \end{equation}
  where $R_n(M, \theta)$ are integrals over $M$ of universal
  polynomials in Tanaka-Webster curvature and covariant derivatives
  (with respect to the classical elliptic development given in
  \cite[Lem 1.7.7]{Gilkey}, the only changes here concern the powers
  of $t$: this is due to the fact that, in the Heisenberg calculus,
  horizontal directions have weight $1$, while $T$ is of weight $2$.
  For instance, this implies that the ``Heisenberg-dimension'' of $M$
  is $4$ instead of $3$).

\begin{rem}    
  Another more direct track, if steeper, also leads to such kernel
  developments. One can follow Beals-Greiner-Stanton's approach to
  heat kernels asymptotics in the contact setting. In \cite{Be-Gr-St}
  they have extended their symbolic calculus on $M × \R$ to
  include the heat operator $\partial_t + P$ for some positive
  sub-Laplacians $P$. They show that in the case $P$ is a positive
  Folland-Stein type operator, one can inverse the symbol of
  $\partial_t + P$ inside this calculus, which gives rather directly
  developments like \eqref{eq:43} for $\Tr(Q\exp^{-tP})$ from the
  symbol of $Q(\partial_t + P)^{-1}$, see also \cite[Sec 4]{Getzler}.
  By R.~Ponge's recent work \cite{Pongeindice,Ponge}, this approach
  leads to a relatively simple proof of the index theorem, and also
  applies to more general positive hypoelliptic $P$ as $(\Delta_2)^2$.
\end{rem}
  
Let us now complete the proof of Theorem \ref{thm:eta(D*)}.  Mellin
transform and the functional calculus relate the asymptotic
development in small time of the heat kernel to $\eta$ and $\zeta$
functions \cite[Section 1.10]{Gilkey}. In particular, \cite[p
81]{Gilkey} and \eqref{eq:43} yield:
$$
\eta(\Delta_2, s) \,\Gamma((s+1)/2) \ = \ \sum_{n=0}^N \frac{4}{2s
  + n - 4} R_n(M, \theta) + h_N(s)
$$
where $h_N$ is an holomorphic function for $s > 2 -N/2$. Hence we
get the required meromorphic extension of $\eta(\Delta_2)(s)$. The
same technique applies to $\Delta_3$ on $\Omega^3 M$, but this is a
positive operator whose heat kernel development has been extensively
treated in \cite[Theorem 7.30]{Be-Gr-St}: the $\eta$-function is here a
$\zeta$-function which is regular at $s=0$.

Using the spectral decomposition \eqref{eq:42}, we get that
$\eta(D*)(s)$ is meromorphic with $s=0$ being possibly a simple pole.
It remains to show that this function is regular at $s=0$. We first
note that the value of the residue of $\eta(D*)$ at $s=0$ is $2
R_4(M,\theta)$. It is easily seen in \eqref{eq:13} that $D*$ becomes
$k\, D*$ in the contact rescaling $\theta \to k \theta$.
Therefore, $\eta (D*_{k\theta})(s) = k^{s} \eta (D*_\theta)(s)$ and
$$
R_4(M, k \theta) = R_4(M, \theta).
$$
Following the proof of Lemma \ref{lem:2:3}, this implies that,
\emph{in dimension $3$},
  \begin{equation}
    \label{eq:44}
    R_4(M, \theta) \ = \ C_1\,\int_M R^2 \, \theta \land d\theta
     \ + \ C_2\,\int_M  |\tau|^2 \,\theta \land d\theta
  \end{equation}
  where $R$ and $\tau$ are Tanaka-Webster curvature and torsion and
  $C_1$, $C_2$ are universal constants.
  
  The residue is moreover invariant under smooth deformation of the
  pseudohermitian and CR structures (\emph{i.e.} both $\theta$ and
  $J$): as underlined in \cite[Lemma 1.10.2]{Gilkey} this general
  feature stems from the existence of a local variation formula for
  $\eta$-functions, namely in the absence of spectral flow here:
  $$
  \dot \eta(\Delta_2)(s) \, = \, -s\, \tr ( \dot\Delta_2
  \Delta_2^{-(s+1)/2}).
  $$
  The point here is that the trace on the right has a meromorphic
  extension coming from the development of $\tr(\dot \Delta_2 \exp^{-t
    (\Delta_2)^2})$, but the possible simple pole at $s=0$ is actually
  cancelled out by the $s$ in front of the whole expression.
  
  The conclusion is that the integrals in \eqref{eq:44} have to be
  independent of variations of $\theta$ and $J$, and this implies
  $C_1= C_2= 0$: indeed, the variations of $R^2$ and $|\tau|^2$ when
  $\theta \to \theta_f= \exp^{2f} \theta$ have been computed
  in \cite[Sec. 5]{Lee86}.  One finds that
  \begin{equation}
    \label{eq:45}
      \frac{d}{df} (R^2 \,\theta \land d\theta ) = 8 R\, ( \Delta_H f
      ) \, \theta \land d\theta
  \end{equation}
  while (if $\tau = A_{11} \theta^1 \otimes\theta^1$)
  \begin{equation}
    \label{eq:46}
      \frac{d}{df} (|\tau|^2\,\theta \land d\theta ) = 2i (A_{\bar
        1\bar 1} f_{,11} - A_{11} f_{,\bar 1 \bar 1})\, \theta \land
      d\theta.  
    \end{equation}
    After integration by parts, this yields
    \begin{equation}
      \label{eq:47}
      \frac{d}{df}R_4(M, \theta) 
      \ = \ 8C_1 \,\int_M f\Delta_H R \, \theta \land d\theta
      \, + \, 2iC_2 \,\int_M f(A_{\bar 1\bar 1,11} - A_{11,\bar 1\bar 1})
      \,\theta \land d\theta.
    \end{equation}
    Testing on a circle bundle (with vanishing torsion) over a Riemann
    surface of non constant curvature cancels out $C_1$. General
    expression for torsion of hypersurfaces in \cite[Sec. 4]{Webster}
    shows that $A_{\bar 1 \bar 1, 11} - A_{11, \bar 1 \bar 1}$ does
    not vanish identically: actually, following \cite{Lee88} the only
    Bianchi identity of order $2$ between $R$ and $\tau$ in dimension
    $3$ is $R_{,0} = A_{11, \bar 1 \bar 1} + A_{\bar 1\bar 1, 11}$,
    which does not occur in \eqref{eq:47} so that $C_2=0$.
\end{proof}

\begin{rem}
  The contact-de Rham complex exists on contact manifolds of any
  dimension, and the contact-signature operator $D*$ is still
  self-adjoint in dimension $4n -1$. Therefore the properties of
  $\eta(D*)(s)$ stated in Theorem \ref{thm:eta(D*)} make sense on
  contact manifolds of any dimension. Most of the previous discussion,
  and its conclusions, still applies, but the last argument about the
  regularity at $s=0$ of $\eta(D*)$.  The residue is still both a
  contact invariant, independent of the choices of $\theta$ and $J$,
  and an integral of some universal pseudohermitian polynomial of the
  right weight. But many possibilities are now left, which cannot be
  so easily analysed (even in the next relevant dimension $7$, the
  algebra becomes quite complicated). At the present time, one still
  ignores whether this residue always vanishes or not.
\end{rem}

\subsection*{The CR invariant correction of $\eta(D*)$}

Having now a well-defined object at hand, we can proceed to the
construction of a \emph{modified contact $\eta$-invariant}.

\begin{theo}\label{maintheo4}
  There exists a unique choice of universal constants $C_1$ and $C_2$
  such that, for any compact strictly pseudoconvex CR
  $3$-manifold $M$, the following pseudohermitian invariant
  \begin{equation}\label{eq:etaB} 
  {\overline\eta} (D*) \ = \ \eta(D*) \, + \, C_1   \,\int_{M} R^2
  \theta\land d\theta \,  
  + \, C_2 \,\int_{M} |\tau|^2 \,\theta\land  
  d\theta \ ,
  \end{equation}
  formed from a contact form $\theta$, its Tanaka-Webster curvature
  $R$ and torsion $\tau$, is in fact a \emph{CR invariant} of $M$,
  which we shall call the \emph{modified contact}
  $\eta$-\emph{invariant}.
\end{theo}

The key point for the proof of Theorem \ref{maintheo4} is the
following: on an oriented CR $3$-manifold $M$, the space of adapted
contact forms for a given CR structure (let us denote it by $\Theta$) is
contractible and non-empty. Then, for a CR invariant, being CR invariant
simply means being independent of the choice of the contact form, \emph{i.e.}
having a vanishing derivative in the direction of any variation in
$\theta$.

Using the analysis above, we get that $\eta(D*)$, seen as a function
on the space $\Theta$ of contact forms adapted to a given CR structure, has the
following features :
\begin{enumerate}
\item $\eta(D*_{k\theta}) = \eta(D*_{\theta})$ for any positive $k$;
\item its derivative is local: if $\theta_t = (1+tf)\theta$ is a small
  variation of contact forms,
\begin{equation*}
  \frac{d}{dt} \eta(D*_{\theta_t}){}_{t=0} = \int_M f \,
  \mathcal{E_{\theta}}    
\,\theta\land d\theta \ ,
\end{equation*}
where $\mathcal{E_{\theta}}$ is a local pseudohermitian invariant of
$\theta$ built algebraically and \emph{universally} from a finite jet
of $\theta$ and its Tanaka-Webster curvature $R$ and torsion $\tau$.
\end{enumerate}
One then deduces from (i) and (ii) that, necessarily,
\begin{equation}\label{eq:invariance-deta1}
\mathcal{E}_{k\theta} = k^{-4} \mathcal{E}_{\theta} \ ,
\end{equation}
and moreover
\begin{equation}\label{eq:invariance-deta2}
\int_M \mathcal{E}_{\theta}\,\theta\land d\theta  \ = \ 0 \ .
\end{equation}
Said otherwise, $\mathcal{E}_{\theta}$ is of weight $-4$ and vanishing
integral.  One can then remark a basic fact:

\begin{lem}\label{lem:derivees}
  Let $\alpha$ be a smooth \emph{closed}, and \emph{real} $1$-form on
  $\Theta$ where $T_{\theta}\Theta$ is identified to the space of
  functions on $M$ through $f \to \frac{d}{dt}(1+tf)\theta$.
  If $\alpha$ is of the type
  \begin{equation}
    \alpha_{\theta} : f \in C^{\infty}(M) \longmapsto \alpha_{\theta}(f) 
    = \int_M f\, \mathcal{A}_{\theta} \, \theta\land d\theta 
  \end{equation}
  where $\mathcal{A}_{\theta}$ is an universal local pseudohermitian
  invariant of a finite jet of $\theta$ of weight $-4$ and vanishing
  integral, then $\alpha$ is a linear combination of the derivatives
  in $\theta$ of
\begin{equation*} 
\int_M R^2\theta\land d\theta \quad \textrm{and} 
\quad \int_M |\tau|^2\theta\land d\theta \ .
\end{equation*}
\end{lem}

\begin{proof}
  We argue as in section \ref{sec:2:3}, classifying local
  pseudo-hermitian invariants that are real and of weight $4$. We have
  seen that the sole possibilities are:
\begin{equation*}
  \begin{split}
    R^2, \quad |\tau|^2,\quad & R_{,0} = A_{11,\bar 1\bar 1} + A_{\bar
      1\bar 1,11}
    \ \ \textrm{(Bianchi identity)} ,\\
    & \Delta_H R,\quad i(A_{11, \bar 1\bar 1} - A_{\bar 1\bar 1,11}) \ 
    .
  \end{split}
\end{equation*}
The first two expressions have non-vanishing integrals in general,
they then have to be forgotten. From \eqref{eq:45} the fourth is the
variation of $\frac{1}{8} \int_M R^2\,\theta\land d\theta$, whereas
from \eqref{eq:46} the fifth is the variation of $\frac{-1}{2} \int_M
|\tau|^2 \,\theta\land d\theta$.

We check that the third one does not yield a closed form. According to
\cite[Sec. 5]{Lee86}, a change of contact form $\theta \to
\theta_f =  e^f \theta$ induces the following changes
$$
R_f = e^{-f} (R + 2 \Delta_H f - 2 |f_{,\bar 1}|^2) \quad
\mathrm{and} \quad T_f = e^{-f} (T + if_1 Z_{\bar 1} - i f_{\bar 1}
Z_1),
$$
and therefore
$$
\frac{d}{d f} (R_{,0} \theta \land d\theta) = \bigl( -f_{,0} R + i
f_{,1} R_{,\bar 1} -i f_{,\bar 1} R_{,1} + 2 (\Delta_H f)_{,0} \bigr) \,
\theta \land d\theta \,.
$$
When restricted on the sphere $\cerc^3$ with its constant
curvature pseudohermitian structure this gives
\begin{align*}
  \int_M (g \frac{d}{df} - f \frac{d}{dg})(R_{,0} \,\theta \land
  d\theta) &= 2\int_M ((\Delta_H f)_{,0} g - (\Delta_H g)_{,0} f )
  \theta \land d\theta \\
  & = -4 \int_M (\Delta_H f) (T.g) \,  \theta \land d\theta.
\end{align*}
This expression does not vanish identically: for instance when taking
any non $T$-invariant function $g$ and $f$ such that $\Delta_H f =
T.g$.  This completes the proof.
\end{proof}
 
This shows Theorem \ref {maintheo4}, exhibiting a new CR invariant
\begin{equation}
  {\overline\eta}(D*) \ = \ \eta(D*)\, + \,C_1 \int_M R^2\, \theta\land
  d\theta  + C_2 \int_M |\tau|^2\,\theta\land d\theta \ .   
\end{equation}
Uniqueness in the choice of the constants is obtained because no
linear combination in the integrals of $R^2$ and $|\tau|^2$ can be a
CR invariant. \qed

\smallskip
\begin{rem}\label{maintheo5}
An analogous line of reasoning yields:
  there exists a universal constant $C'$ such that, for any compact
  strictly pseudoconvex Cauchy-Riemann $3$-manifold $M$,
  \begin{equation}\label{eq:nuetaB} 
  {\overline\eta}(D*) \ -\ C'\, \nu(M)
  \end{equation}
  is a \emph{contact} invariant, \emph{i.e.} is independent of the
  choice of the complex structure. The proof (left to the reader)
  consists in proving that the only tensorial choice for the
  differential of $\overline{\eta}$ is (up to some multiplicative
  constant) the Cartan curvature like in \eqref{eq:3} and
  \eqref{eq:4}.

  Of course, in view of the relation (\ref{eq:etanugeneral}) between
  $\nu$ and $\eta(D*)$ in the CR-Seifert case, one expects that the
  constants $C'$ above and and $C_1$ in Theorem \ref{maintheo4} should
  be respectively $-\frac{1}{3}$ and
  $(\frac{1}{512}-\frac{1}{48\pi^2})$, but the case of CR-Seifert
  manifolds is not sufficient to determine them.
  The best one can get is the following: it has already been remarked
  earlier that the value of the renormalized $\eta$-invariant $\eta_0$ is
  purely topological on CR-Seifert manifolds. Keeping the contact form
  fixed, this means that it has to be independent of the complex
  structure.  As $\eta(D*) = \eta_0 - \frac{1}{512}\int R^2\theta\land d\theta$ and
$$
{\overline\eta} - C'\nu = (1+3C')\eta_0 + (C_1 - \frac{1}{512} - \frac{C'}{16\pi^2}) 
\int R^2 \theta\land d\theta
$$ 
must be a contact invariant, this implies that
$$ 
C_1 - \frac{1}{512} - \frac{C'}{16\pi^2} = 0 ,
$$
since the integral of $R^2$ has non-zero variations with respect to the
complex structures. 

Guessing the values of $C$ in Conjecture
\ref{mainconj} and $C_2$ in Theorem \ref{maintheo4} seems much harder.
Having a precise value for them would (for instance) involve a precise
computation of the spectrum of $\eta(D*)$ in a case where the torsion
does not vanish. This seems difficult to achieve either with our
methods, which rely on Fourier decomposition under the circle action,
or with classical tools of representation theory, which require a high
degree of homogeneity.

Of course, one knows that the derivative of $\eta(D*)$ is given by algebraic 
expressions of the jet of the hypoelliptic symbols of the involved operators. 
However these expressions are so intricate that the constants are only computable 
this way ``in theory'', and not in practice.
\end{rem}

\begin{rem}
  The same arguments also apply to the renormalized $\eta$-invariant
  $\eta_0$ introduced in section \ref{sec:adia}, instead
  of $\eta(D*)$. This explains \emph{a priori} the existence of some local
  correction of $\eta_0$ leading to a CR invariant, itself related (up
  to some contact invariant) to a multiple of $\nu$; this might be compared
  with Lemma \ref{lem:2:3}.
\end{rem}

\bigskip

\section{Proof of the corollaries}\label{sec:coro}

Corollaries \ref{cor:1} and \ref{cor:3} rely on the formula discovered
by the first and second authors \cite[Theorem 1.2]{ob-mh1}: for any
Einstein asymptotically hyperbolic manifold $(N^4,g)$,
\begin{equation} \label{eq:my}
\frac{1}{8\pi^2} \int_N \left( 3|W^-|^2 - |W^+|^2 + \frac{1}{24}
\Scal^2 \right) - \chi(N) + 3\,\tau(N) = \nu(M).
\end{equation}
For complex hyperbolic surfaces, the integral term is zero.  If $\bar
N$ is smooth, with $M$ as the only end, then the topological
contributions always are integers.  Corollary \ref{cor:1} is then
proved.

\medskip

It is instructive to check the results for a holomorphic disk bundle
over a hyperbolic Riemann surface $\Sigma$, with $M$ as its boundary.
Clearly one has $\chi(N) = \chi(\Sigma)= \chi$ and $\tau(N) =-1$.  If
$N$ carries a complex hyperbolic metric with $M$ as its boundary at
infinity, then corollary \ref{cor:1} gives the equation
$$
\chi\, - \, 3\,\tau \ = \ -\,\nu(M) \ = \ d + 3 + \frac{\chi^2}{4d}
$$
and the only solution is $d=\frac{\chi}{2}$.  We then recover the
well-known fact that the only disk bundles carrying a complex
hyperbolic metric are the square roots of the (complex) tangent
bundle.

\medskip

Corollary \ref{cor:3} is again a direct consequence of (\ref{eq:my}),
since for a Kähler-Einstein metric, the integral term is non negative.
For an Einstein metric, the story is more complicated, but positivity
is achieved if solutions of the Seiberg-Witten equations exist, and it
is proved in \cite[corollary 31]{rollin} that it is a consequence of
the nonvanishing of the Kronheimer-Mrowka invariants \cite{kro-mro}.

\medskip

From \cite[Theorem 5.12]{carron-pedon}, one knows that pseudoconvex
complex hyperbolic surfaces $N$ have vanishing third homology group
$H_3(N,\mathbb{Z})$. Hence no multiple ends can occur, but one expects
orbifold singularities or cusps to appear in the interior of a complex
hyperbolic filling. The complex hyperbolic cusps can be compactified
to yield a complex orbifold surface that we note again $N$, by adding
at the infinity of each cusp a quotient $\Sigma_i$ of a 2-torus.  The
Corollaries \ref{cor:1} and \ref{cor:3} remain true in this case, with
the Euler characteristic and the signature of $N$ being replaced by
their orbifold versions: In case $\ell$ cusps are present, there is an
additional contribution in the signature coming from the
self-intersection of each 2-torus at infinity.  Namely, one has to
consider the modified signature \cite[proposition
3.4]{biquard-rigidite-fini}
$$
\tau_{\mathrm{cusp}}(N)=\tau(N) - \frac{1}{3}\sum_1^\ell
[\Sigma_i]\cdot[\Sigma_i] .
$$
Of course, Corollary \ref{cor:2} is no more true, since the
characteristic numbers are now rational; the denominator of $\nu$ only
gives an hint on the order of the singularities needed to fill $M$.

\medskip

\subsection*{Explicitation for lens spaces} 
We now specialize the formula obtained in Corollary \ref{maintheo3} to
the lens space $L(p,q)$ obtained as a quotient of the 3-sphere
$\cerc^3$ in $\mathbb{C}^2$ by $\setZ/p\setZ$, with its generator
acting on $\mathbb{C}^2$ by
$(\exp^{\frac{2i\pi}{p}},\exp^{\frac{2iq\pi}{p}})$, where $q$ is prime
with $p$. They are interesting in connection with filling by Einstein
metrics, since some of them appear as boundary at infinity of selfdual
Einstein metrics \cite{cal-sin}. On the other hand, it has been shown
that large families of them admit symplectic fillings \cite{lisca}, so
that Corollary \ref{cor:3} may be applied to these.

\begin{prop}\label{nu:lens} One has:
  $\nu(L(p,q))=-\frac{1}{p}+12\, s(p,q,1)$.
\end{prop}

For sake of comparison, we recall to the interested reader the value
of the classical $\eta$-invariant on lens spaces with the standard
round metric, as computed by Atiyah-Patodi-Singer \cite[Proposition
2.12]{aps2}:
\begin{equation}\label{eq:aps}
\eta(L(p,q)) = -4\, s(p,q,1) .
\end{equation}

\begin{proof}
  For simplicity, we shall assume that $(q-1)$ is prime with $p$ (as a
  matter of fact this implies that we take $q\neq 1$), and we leave
  the general case to the reader.  Let us see the 3-sphere as the
  bundle $\mathcal{O}(-1)$ over the projective line $\mathbb{C}P^1$.
  The induced action on $\mathbb{C}P^1$ has two fixed points: the two
  antipodal points, with action of $\setZ/p\setZ$ generated by
  $\exp^{± i2\pi\frac{q-1}{p}}$, and action in the fiber by
  $\exp^{i\frac{2\pi}{p}}$ and $\exp^{i2\pi\frac{q}{p}}$ respectively.
  Therefore $L(p,q)$ is a $\cerc^1$-orbifold bundle over an orbifold
  projective line with two orbifold points with angle
  $\frac{2\pi}{p}$. The Euler characteristic is $\chi=\frac{2}{p}$ and
  the degree (first Chern number) is $d=-\frac{1}{p}$. Now Corollary 
  \ref{maintheo3} and Ouyang's Theorem \ref{th:ouy} give the formulae
\begin{align*}
  \nu(L(p,q))&=-3+\frac{2}{p}-12\big(s(p,q-1,1)+s(p,1-q,q)\big)\ , \\
  \eta(L(p,q))&=1-\frac{1}{p}+4\big(s(p,q-1,1)+s(p,1-q,q)\big)\ ,
\end{align*}
(note that the extra parameter $\rho$ in Theorem \ref{th:ouy}
appears naturally on lens spaces), 
so that $\nu(L(p,q))=-\frac{1}{p}-3\eta(L(p,q))$. The proposition then
follows from (\ref{eq:aps}).
\end{proof}

\smallskip

\subsection*{Comparison with the Burns-Epstein invariant}
Another interesting point is to compare these results with those
obtained by use of the Burns-Epstein $\mu$-invariant \cite{BE88,BE90b}
(it is already suggested at the end of \cite{BE90b} that obstructions
follow from computations of $\mu$).  The $\mu$-invariant is defined on
strictly pseudoconvex CR 3-manifolds with trivial tangent holomorphic
bundle only. Roughly speaking, it comes from Chern-Simons-type
constructions (integration of a local formula), whereas the
$\nu$-invariant is extracted from the Atiyah-Patodi-Singer
$\eta$-invariant.  The relation between $\mu$ and $\nu$ is similar to
that between the $\eta$ and the Chern-Simons invariants: more
precisely, when $\mu$ is defined, then for a CR structure $J$ one has
$$
\nu(J) \ = \ 3\,\mu(J) \ + \ \text{constant} , 
$$
with the constant depending only of the underlying contact
structure \cite[Theorem 1.3]{ob-mh1}.  Burns-Epstein's version of
Miyaoka-Yau \cite{BE90b} then reads, if $M$ is the boundary at
infinity of a Kähler-Einstein $N$:
\begin{equation}\label{eq:BE}
\chi(N)-\frac{1}{3}\,\bar{c}_1(N)^2 \geqslant - \mu(M) ,
\end{equation}
with equality if the metric is complex hyperbolic; here $\bar{c}_1$ is
a lift in $H^2(N,M)$ of $c_1(N)$.

A first important difference here is that our obstruction in Corollary
\ref{cor:3} (filling by an ACH Einstein metric) is purely topological,
whereas (\ref{eq:BE}) involves a complex structure and a
Kähler-Einstein metric.

Another important fact to be noticed, at least in the case when the
quotient has no orbifold singularities, is that the obstructions
obtained by both methods are different: if $M$ is a $\cerc^1$-bundle
over the Riemann surface $\Sigma$, then the $\mu$-invariant, being
defined by a local formula, is multiplicative on finite coverings
\cite{BE88,BE90b}.  Hence the values are
\begin{equation}\label{eq:mu.vs.nu} 
\mu = \frac{\chi^2}{4d} 
\quad \textrm{whereas} \quad  \nu = - \frac{ \chi^2}{4d} - d - 3.
\end{equation}
Equation (\ref{eq:BE}) implies that $3\mu$ must be an integer, {\it
  i.e.}  $\frac{3\chi^2}{4d}$ must belong to $\mathbb{Z}$, a condition
that is weaker than Corollary \ref{cor:2}, by a factor $3$.

\vskip 1cm

\begin{small}
  {\flushleft\sl Acknowledgements}. The authors are grateful to
  Yoshinobu Kamishima for useful conversations on the possible
  applications of $\nu$, and to Elisha Falbel for comments. M.~H. thanks
  Emmanuel Royer for his help in computations of examples at a very
  early stage of this paper.  Finally, we thank Nigel Hitchin for
  inventing the word `diabatic'.
\end{small}

\vskip 1cm

\bibliographystyle{smfplain}

\end{document}